\documentclass[11pt]{article}

\usepackage{fullwidth}
\usepackage[top=1.5in, bottom=1.5in, left=1in, right=1in]{geometry}
\usepackage{xypic}
\usepackage{amsgen, amsmath, amstext, amsbsy, amsopn, amsfonts, amssymb, comment}
\usepackage{eepic}
\usepackage{epsf}
\xyoption{all}
\usepackage{tikz}
\usetikzlibrary{calc,through,backgrounds,shapes,matrix}

\usepackage[backref=page]{hyperref}
\usepackage{url, hypcap}
\hypersetup{colorlinks=true, citecolor=darkblue, linkcolor=darkblue}
\usepackage{mathrsfs}

\definecolor{darkblue}{rgb}{0.0,0,0.7} 
\definecolor{darkred}{rgb}{0.7,0,0} 
\newcommand{\darkred}{\color{darkred}} 
\definecolor{lightgrey}{rgb}{0.7,0.7,0.7} 

\newcommand{\defn}[1]{\emph{\darkred #1}} 

\def\ds{\displaystyle}
\def\Box{\square}

\def\mapright#1{\smash{\mathop{\longrightarrow}\limits^{#1}}}
\def\tra#1{\smash{\mathop{\mid\kern
-1pt\joinrel\relbar\joinrel\relbar}\limits^{*}_{#1}}}
\def\longtra#1{\smash{\mathop{\mid\kern
-1pt\joinrel\relbar\joinrel\relbar\joinrel\relbar}\limits^{*}_{#1}}}
\def\vlongtra#1{\smash{\mathop{\mid\kern
-1pt\joinrel\relbar\joinrel\relbar\joinrel\relbar\joinrel\relbar}\limits^{*}_{#1}}}
\def\vvlongtra#1{\smash{\mathop{\mid\kern
-1pt\joinrel\relbar\joinrel\relbar\joinrel\relbar\joinrel\relbar\joinrel\relbar}\limits^{*}_{#1}}}
\def\vvvlongtra#1{\smash{\mathop{\mid\kern
-1pt\joinrel\relbar\joinrel\relbar\joinrel\relbar\joinrel\relbar\joinrel\relbar\joinrel\relbar}\limits^{*}_{#1}}}
\def\etra#1{\smash{\mathop{\mid\kern
-1pt\joinrel\relbar\joinrel\relbar}\limits_{#1}}}

\def\vlongrightarrow{\relbar\joinrel\longrightarrow}

\def\vvlongrightarrow{\relbar\joinrel\vlongrightarrow}

\def\longmapright#1{\smash{\mathop{\vlongrightarrow}\limits^{#1}}}
\def\vlongmapright#1{\smash{\mathop{\vvlongrightarrow}\limits^{#1}}}

\def\iff{\Leftrightarrow}
\def\Rw{\Rightarrow}
\def\oo{\overline}

\def\wh{\widehat}

\def\I{{\cal{I}}}

\def\N{\mathbb{N}}

\def\res{\operatorname{Res}}
\def\rg{\operatorname{RG}}

\def\rc{\operatorname{RC}}
\def\crc{\operatorname{CRC}}
\def\orc{\operatorname{ORC}}
\def\prc{\operatorname{PRC}}
\def\src{\operatorname{SRC}}

\def\lcs{\mbox{lcs}}

\def\ns{\mbox{NS}}

\def\max{\mbox{max}}
\def\min{\mbox{min}}

\def\tape{\mbox{tape}}
\def\heads{\mbox{heads}}
\def\leg{\mbox{Leg}}
\def\rres{\mbox{RRes}}
\def\lres{\mbox{LRes}}
\def\RSC{\mbox{RSC}}
\def\LSC{\mbox{LSC}}

\def\cay{\mbox{Cay}}

\def\D{{\bf D}}
\def\V{{\bf V}}

\def\p{\varphi}

\def\inv{^{-1}}

\def\bi{\begin{itemize}}
\def\ei{\end{itemize}}
\def\beq{\begin{equation}}
\def\eeq{\end{equation}}

\newtheorem{T}{Theorem}[section]
\newcommand{\bt}{\begin{T}}
\newcommand{\et}{\end{T}}
\newcommand{\ftd}{$\square$\end{T}}

\newtheorem{Proposition}[T]{Proposition}
\newcommand{\bp}{\begin{Proposition}}
\newcommand{\ep}{\end{Proposition}}
\newcommand{\fpd}{$\square$\end{Proposition}}

\newtheorem{Lemma}[T]{Lemma}
\newcommand{\bl}{\begin{Lemma}}
\newcommand{\el}{\end{Lemma}}
\newcommand{\fld}{$\square$\end{Lemma}}

\newtheorem{Corol}[T]{Corollary}
\newcommand{\bc}{\begin{Corol}}
\newcommand{\ec}{\end{Corol}}
\newcommand{\fcd}{$\square$\end{Corol}}

\newtheorem{Remark}[T]{Remark}
\newcommand{\br}{\begin{Remark}}
\newcommand{\er}{\end{Remark}}
\newcommand{\frd}{$\square$\end{Remark}}

\newtheorem{Example}[T]{Example}
\newcommand{\be}{\begin{Example}}
\newcommand{\ee}{\end{Example}}

\newtheorem{Problem}[T]{Problem}
\newcommand{\bq}{\begin{Problem}}
\newcommand{\eq}{\end{Problem}}

\newcommand{\proof}
   {\par\medbreak\noindent{\bf Proof}.\enspace}

\newcommand{\qed}{
$\Box$
\par\bigbreak}

\normalbaselines
\def\abstract#1{\par\bigskip
\begingroup\small
\baselineskip=12truept
\begin{center}ABSTRACT\end{center}
\par\medskip\par\noindent
\null\hfill\hbox{\vbox{\hsize=5truein\noindent#1}}
\hfill\null\par\endgroup\par}

\numberwithin{figure}{section}
\numberwithin{equation}{section}

\title{The semaphore codes attached to a Turing machine via resets and their various limits}
\author{{\bf John Rhodes}\\ 
{\em Department of Mathematics, University of California, Berkeley,}\\ 
{\em CA 94720, U.S.A.}\\
{\em email:} rhodes@math.berkeley.edu, BlvdBastille@aol.com\\
$ $\\
{\bf Anne Schilling}\\
{\em Department of Mathematics, University of California, Davis,}\\
{\em One Shields Ave., Davis, CA 95616-8633, U.S.A.}\\
{\em email:} anne@math.ucdavis.edu\\
$ $\\
{\bf Pedro V. Silva}\\
{\em Centro de
Matem\'{a}tica, Faculdade de Ci\^{e}ncias, Universidade do
Porto,}\\ {\em R. Campo Alegre 687, 4169-007 Porto, Portugal}\\
{\em email:} pvsilva@fc.up.pt}

\date{\today}

\begin{document}

\maketitle

\begin{center}\small
\end{center}

\abstract{
We introduce semaphore codes associated to a Turing machine via resets. 
Semaphore codes provide an approximation theory for resets. In this paper we
generalize the set-up of our previous paper ``Random walks on semaphore codes 
and delay de Bruijn semigroups'' to the infinite case by taking the profinite limit
of $k$-resets to obtain $(-\omega)$-resets. We mention how this opens new
avenues to attack the P versus NP problem.
}

\section{Introduction}

In our previous paper~\cite{RSS:2016}, we developed algebraic foundations centered around the prime 
decomposition theory for finite semigroups and finite automata (see \cite{KR:1965},~\cite{Rhodes:2010}
and~\cite[Chapter 4]{RS:2009}). This analysis focused on the right zero component action, when the corresponding 
pseudovariety contains all finite transformation semigroups $(X,S)$ (or automata) with the property that there exists 
some $k \geq 1$ such that every product $s_1\cdots s_k$ (with $s_i \in S$) is a constant map (or reset) on $X$. 
Such finite automata were called $k$-reset graphs in~\cite[Section 3]{RSS:2016} and their elementary properties 
were studied, using the lattice of right congruences on the finite free objects (De Bruijn semigroups), 
in~\cite[Sections 5 and~6]{RSS:2016}. 

Semaphore codes, which are well known in the literature (see~\cite{BPR:2010}), were 
proved~\cite[Sections 4 and 7]{RSS:2016} to be in bijection with special right congruences and provide a lower 
approximation to any right congruence with the same hitting time to constant. Thus in many applications right 
congruences can be replaced by semaphore codes. 
Except for Section~4 on semaphore codes, all the material in~\cite{RSS:2016} up to and including Section~7 
is restricted to finite codes and automata. Finally, in~\cite[Section 8]{RSS:2016}, a natural random walk on any 
(finite or infinite) semaphore code was constructed and its stationary distribution plus hitting time to constant 
were computed. 

In this paper we do the following: first we reveal a main application we have in mind~\cite{RS:2008} by introducing
the infinite and finite semaphore codes associated to a Turing machine via resets (see Section~\ref{section.TM}). 
Then Sections~\ref{section.free profinite} and~\ref{section.regr2} take the profinite limits of $k$-reset graphs 
yielding $(-\omega)$-reset graphs. 

We consider the pseudovariety {\bf D} and left infinite words, but by duality we have analogous results for the 
pseudovariety {\bf K} and right infinite words. We need both versions for studying Turing machines.
Generalizing the finite case from~\cite{RSS:2016}, we study right congruences and special right congruences 
in bijection with infinite semaphore codes and the natural action in Sections~\ref{section.rc omega} 
and~\ref{section.src omega} and obtain an approximation theory as in the finite case of $k$-resets, but for 
$(-\omega)$-resets. In the final Section~\ref{section.future}, we make some more remarks on relating 
Section~\ref{section.TM} to Sections~\ref{section.free profinite}-\ref{section.src omega}, and the next paper on 
attacking P versus NP.

\subsection*{Acknowledgements}

We are indebted to Jean-Camille Birget for good advice and patient reading. We thank also Benjamin Steinberg 
and Nicolas M. Thi\'ery for discussions.

The first author thanks the Simons Foundation--Collaboration Grants for Mathematicians for  travel grant $\#313548$.
The second author was partially supported by NSF grants OCI--1147247 and DMS--1500050.
The third author was partially supported by CMUP (UID/MAT/00144/2013), which is funded by FCT (Portugal) with 
national (MEC) and European structural funds (FEDER), under the partnership agreement PT2020.

\section{Resets and Turing machines}
\label{section.TM}

In this section we present a new viewpoint on Turing machines centered in the concept of resets and their associated semaphore codes.

\subsection{Turing machines}
\label{geco}

For details on Turing machines, the reader is referred to \cite{HU:1979}.

Let us assume that $T = (Q,A,\Gamma, q_0,F,\delta)$ is a (deterministic) Turing machine, where:
\bi
\item
$Q$ is the (finite) state set;
\item
$A$ is the (finite) input alphabet;
\item
$\Gamma$ is the (finite) tape alphabet (containing $A$ and the blank symbol $B$);
\item
$q_0 \in Q$ is the initial state;
\item
$F \subseteq Q$ is the set of final states;
\item
$\delta: Q \times \Gamma \to Q \times (\Gamma \setminus \{ B \}) \times \{ L,R\}$ is the (partial) transition function.
\ei
Then we write $\Omega = \Gamma \cup (\Gamma \times Q)$. To make notation lighter, we shall denote 
$(X,q) \in \Gamma \times Q$ by $X^q$. To avoid confusion with powers of $X$, we stipulate that from now on 
symbols such as $q,q',q_i$ will be reserved to denote states and never integers.

We define two homomorphisms $\tape\colon \Omega^* \to \Gamma^*$ and $\heads\colon\Omega^* \to (\mathbb{N},+)$ by
$$\tape(X) = \tape (X^q) = X,\quad \heads(X) = 0,\quad \heads(X^q) = 1$$
for all $X \in \Gamma$ and $q \in Q$. Now we define the set of all \defn{legal words} by
$$\leg(T) = B^*\{ w \in \Omega^* \mid \tape(w) \in \{1,B\}(\Gamma\setminus \{ B \})^*\{1,B\}, \; \heads(w) \leq 1\}B^*.$$
Note that $\leg(T)$ is closed under reversal and factors.
The \defn{illegal words} are the elements of the complement $\Omega^* \setminus \leg(T)$, which is the ideal 
having as generating set all the words of the form
\bi
\item
$X_1^quX_2^{q'}$
\item
$YB^nY'$
\item
$YB^mB^qB^kY'$
\item
$Y^qB^nY'$
\item
$Y'B^nY^{q}$
\item
$B^qB^nY$
\item
$YB^nB^q$
\ei
where $X_1,X_2 \in \Gamma$; $Y,Y' \in \Gamma \setminus \{ B\}$; $u \in \Gamma^*$; $n \geq 1$; $m,k \geq 0$.

Note that legal words do not correspond necessarily to the possible content of the tape during a computation 
(or a factor of that content), but they contain such words as particular cases.

We define the \defn{one-move mapping} $\beta\colon \leg(T) \to \leg(T)$ as follows. Given $w \in \leg(T)$, then $\beta(w)$ 
is intended to be obtained from $w$ by performing one single move of $T$ on a tape with content $w$; if $T$ admits 
no such move from $w$ (in particular, if $\heads(w) = 0$), we set $\beta(w) = w$. In all cases except (A) and (B) below, 
the interpretation of $\beta(w)$ is clear and $|\beta(w)| = |w|$. The following two cases deserve extra clarification:
\bi
\item[(A)]
if $w = X^qw'$ and $\delta(q,X) = (\ldots,\ldots,L)$;
\item[(B)]
if $w = w'X^q$ and $\delta(q,X) = (\ldots,\ldots,R)$.
\ei
In these cases, we interpret $\beta(w)$ as $\beta(Bw)$ (respectively $\beta(wB)$), falling into the general case. 
We say that legal words of types (A) and (B) are $\beta$-\defn{singular}. Note that $|\beta(w)| = |w|+1$ in the 
$\beta$-singular case. Note that, by padding the input sequences with sufficiently many $B$'s before and after, 
cases (A) and (B) never occur.
 
For every $n \geq 0$, we denote by $\beta^n$ the $n$-fold composition if $\beta$. We define also a partial 
mapping 
$$\beta^{(n)}\colon \Omega^* \times \Omega \times \Omega^* \to \Omega$$
as follows. 

Let $u,v \in \Omega^*$ and $X \in \Omega$. If $uXv \in \leg(T)$, then $\beta^{(n)}(u,X,v)$ is the symbol replacing 
$X$ at the designated position in the tape after applying $\beta$ $n$ times. If $uXv \notin \leg(T)$, then 
$\beta^{(n)}(u,X,v)$ is undefined.

We say that $T$ is \defn{legal-halting} if, for every $u \in \leg(T)$, the sequence $(\beta^n(u))_n$ is eventually 
constant. This implies that the sequence $(\beta^{(n)}(u,X,v))_n$ is also eventually constant for all $u,v \in \Omega^*$ 
and $X \in \Omega$ such that $uXv \in \leg(T)$. We write
$$\beta^{\omega}(u) = \lim_{n\to\infty}\beta^n(u), \quad \beta^{(\omega)}(u,X,v) = \lim_{n\to\infty}\beta^{(n)}(u,X,v).$$
Note that, from a formal viewpoint, a Turing machine which halts for every input is not necessarily legal-halting (since 
not every legal word arises from an input configuration), but it can be made legal-halting with minimum adaptations. 

\subsection{Resets}

We say that $r \in \Omega^*$ is a \defn{right reset} if 
$$\beta^{(n)}(t_1rt_2,X,t_3) = \beta^{(n)}(t'_1rt_2,X,t_3)$$
for all $t_1,t'_1,t_2,t_3 \in \Omega^*$, $X \in \Omega$ and $n \geq 0$ such that both 
$t_1rt_2Xt_3, t'_1rt_2Xt_3 \in \leg(T)$. Let $\rres(T)$ denote the set of all right resets of $T$. 

Dually, $r \in \Omega^*$ is a \defn{left reset} if 
$$\beta^{(n)}(t_1,X,t_2rt_3) = \beta^{(n)}(t_1,X,t_2rt'_3)$$
for all $t_1,t_2,t_3,t'_3 \in \Omega^*$, $X \in \Omega$ and $n \geq 0$ such that both 
$t_1Xt_2rt_3, t_1Xt_2rt'_3 \in \leg(T)$. Let $\lres(T)$ denote the set of all left resets of $T$. 

\bl
\label{ideal}
{\rm RRes}$(T)$ and {\rm LRes}$(T)$ are ideals of $\Omega^*$ containing all the illegal words.
\el

\proof
We prove the claim for right resets.

If $r$ is illegal then $t_1rt_2Xt_3$ is always illegal and so $r \in \rres(T)$ trivially. In particular, $\rres(T)$ is nonempty.

Let $r \in \rres(T)$ and let $x,y \in \Omega^*$. Suppose that $xry \notin \rres(T)$. Then there exist 
$t_1,t'_1,t_2,t_3 \in \Omega^*$, $X \in \Omega$ and $n \geq 0$ such that both $t_1xryt_2Xt_3$ and 
$t'_1xryt_2Xt_3$ are legal and
$$\beta^{(n)}(t_1(xry)t_2,X,t_3) \neq \beta^{(n)}(t'_1(xry)t_2,X,t_3).$$
Rewriting this inequality as
$$\beta^{(n)}((t_1x)r(yt_2),X,t_3) \neq \beta^{(n)}((t'_1x)r(yt_2),X,t_3),$$
we deduce that $r \notin \rres(T)$, a contradiction. Thus $xry \in \rres(T)$ and so $\rres(T)$ is an ideal of $\Omega^*$.
\qed

\bl
\label{brl}
If $u \in \Omega^* \setminus B^*$, then $Bu \in {\rm RRes}(T)$ and $uB \in {\rm LRes}(T)$.
\el

\proof
It is easy to see that $Bu$ is a right reset since the only legal words of the form $tBut'$ must arise from $t \in B^*$. 
Similarly, $uB$ is a left reset.
\qed

\bl
\label{betres}\mbox{}
\bi
\item[(i)] $\beta({\rm RRes}(T)) \subseteq {\rm RRes}(T)$.
\item[(ii)] $\beta({\rm LRes}(T)) \subseteq {\rm LRes}(T)$.
\ei
\el

\proof
We prove the claim for right resets.

Let $r \in \rres(T)$. 
We may assume that $r$ is legal, $X^q$ occurs in $r$ and $\delta(X,q)$ is defined. Let $t_1,t'_1,t_2,t_3 \in \Omega^*$, 
$X \in \Omega$ and $n \geq 0$ be such that both $t_1\beta(r)t_2Xt_3$ and $t'_1\beta(r)t_2Xt_3$ are legal. 

Suppose first that $r$ is not $\beta$-singular. Then $t_1rt_2Xt_3$ and $t'_1rt_2Xt_3$ are also legal, and
$$\beta^{(n)}(t_1\beta(r)t_2,X,t_3) = \beta^{(n+1)}(t_1rt_2,X,t_3) = \beta^{(n+1)}(t'_1rt_2,X,t_3) = \beta^{(n)}(t'_1\beta(r)t_2,X,t_3).$$

Thus we may assume that $r$ is $\beta$-singular. Suppose first that
$r = X^qr'$ and $\delta(q,X) = (p,Y,L)$. Then $\beta(r) = B^pYr'$. Since 
$t_1B^pYr't_2Xt_3$ and $t'_1B^pYr't_2Xt_3$ are legal, it is easy to check that $t_1Brt_2Xt_3$ and $t'_1Brt_2Xt_3$ 
are legal as well. Hence
$$\beta^{(n)}(t_1\beta(r)t_2,X,t_3) = \beta^{(n+1)}(t_1Brt_2,X,t_3) = \beta^{(n+1)}(t'_1Brt_2,X,t_3) = \beta^{(n)}(t'_1\beta(r)t_2,X,t_3).$$

Finally, suppose that
$r = r'X^q$ and $\delta(q,X) = (p,Y,R)$. Then $\beta(r) = r'YB^p$. Since 
$t_1r'YB^pt_2Xt_3$ and $t'_1r'YB^pt_2Xt_3$ are legal, it is easy to check that $t_1rBt_2Xt_3$ and $t'_1rBt_2Xt_3$ 
are legal as well. Hence
$$\beta^{(n)}(t_1\beta(r)t_2,X,t_3) = \beta^{(n+1)}(t_1rBt_2,X,t_3) = \beta^{(n+1)}(t'_1rBt_2,X,t_3) = \beta^{(n)}(t'_1\beta(r)t_2,X,t_3).$$

Therefore $\beta(r)$ is a right reset in any case.
\qed

\subsection{Semaphore codes and the output function}

Given an alphabet $X$, we define the \defn{suffix order} on $X^*$ by
$$\mbox{$u \leq _s v$ if $v \in X^*u$}.$$ 
Dually, we define the \defn{prefix order} $\leq_p$.

We say that $S \subseteq X^*$ is a \defn{(right) semaphore code} if $S$ is a suffix code (i.e. an antichain for the suffix 
order) and $SX \subseteq X^*S$. By~\cite[Proposition 4.3]{RSS:2016}, $S$ is a semaphore code if and only if $S$ is 
the set of $\leq_s$-minimal elements in some ideal $I \unlhd X^*$, denoted by $I\beta_{\ell}$.  

Dually, $S \subseteq X^*$ is a \defn{left semaphore code} if $S$ is a prefix code (i.e. an antichain for the prefix order) 
and $XS \subseteq SX^*$. Then $S$ is a left semaphore code if and only if $S$ is the set of $\leq_p$-minimal elements 
in some ideal of $X^*$.

We describe now the semaphore code $\RSC(T)$ defined by the right resets. It
consists of all minimal right resets for the suffix order. If $1 \notin \rres(T)$ (a trivial case), then $\RSC(T)$ consists of 
all right resets $Xz$ with $X \in \Omega$ and $z \in \Omega^*$ such that $z \notin \rres(T)$. In particular, $z$ must 
be a legal word. 

Similarly, the left semaphore code $\LSC(T)$ consists of all minimal left resets for the prefix order. If 
$1 \notin \lres(T)$, then $\LSC(T)$ consists  of all left resets $zX$ with $z \in \Omega^*$ and $X \in \Omega$ 
such that $z \notin \lres(T)$. In particular, $z$ must be a legal word. 

Note also that every right reset contains some $s \in \RSC(T)$ as a suffix. Dually, every left reset contains some 
$s \in \LSC(T)$ as a prefix. 

From now on, we assume that $T$ is legal-halting. The \defn{output function} $\p_T$ is the restriction of the partial 
function $\beta^{(\omega)}:\Omega^* \times \Omega \times \Omega^*$ to 
$(\RSC(T) \cup \{ 1 \}) \times \Omega \times (\LSC(T) \cup \{ 1 \})$.

\bp
\label{ofdt}
Let $T$ be a legal-halting Turing machine. Then $\beta^{\omega}$ is fully determined by the output function $\p_T$.
\ep

\proof
Clearly, $\beta^{\omega}$ is fully determined by $\beta^{(\omega)}$. Let $u,v \in \Omega^*$ and $X \in \Omega$. 
If $uXv$ is illegal, then $\beta^{(\omega)}(u,X,v)$ is undefined, hence we may assume that $uXv \in \leg(T)$. It follows 
that also $BuXvB \in \leg(T)$ and $\beta^{(\omega)}(Bu,X,vB) = \beta^{(\omega)}(u,X,v)$. 

If $u \in B^*$, then $\beta^{(\omega)}(Bu,X,vB) = \beta^{(\omega)}(1,X,vB)$. If $u \notin B^*$, then $Bu \in \rres(T)$ 
by Lemma \ref{brl} and so $Bu = xr$ for some $x \in \Omega^*$ and $r \in \RSC(T)$. It follows that 
$\beta^{(\omega)}(Bu,X,vB) = \beta^{(\omega)}(r,X,vB)$, so in any case we have 
\beq
\label{ofdt1}
\beta^{(\omega)}(Bu,X,vB) = \beta^{(\omega)}(r,X,vB)
\eeq
for some $r \in \RSC(T) \cup \{ 1 \}$.

Now if $v \in B^*$, then $\beta^{(\omega)}(r,X,vB) = \beta^{(\omega)}(r,X,1) = \p_T(r,X,1)$, hence we may assume 
that $v \notin B^*$. Then $vB \in \lres(T)$ by Lemma \ref{brl} and so $vB = r'x'$ for some $r' \in \LSC(T)$ and 
$x \in \Omega^*$. It follows that $\beta^{(\omega)}(r,X,vB) = \beta^{(\omega)}(r,X,r') = \p_T(r,X,r')$, so in view 
of~\eqref{ofdt1} we have that $\beta^{(\omega)}(u,X,v) = \beta^{(\omega)}(Bu,X,vB)$ is determined by $\p_T$. 
Therefore $\beta^{\omega}$ is determined by $\p_T$.
\qed

\subsection{Length restrictions}

For every $\ell \geq 0$, we define the cofinite ideals 
$$\begin{array}{l}
\rres_{\ell}(T) = \rres(T) \cup \Omega^{\ell}\Omega^*,\\
\lres_{\ell}(T) = \lres(T) \cup \Omega^{\ell}\Omega^*.
\end{array}$$
Note that
\beq
\label{prima1}
\displaystyle
\rres(T) = \bigcap_{\ell \geq 0} \rres_{\ell}(T), \quad
\lres(T) = \bigcap_{\ell \geq 0} \lres_{\ell}(T).
\eeq

Since 
$$\beta(\Omega^{\ell}) \subseteq \Omega^{\ell} \cup \Omega^{\ell+1},$$
it follows from Lemma \ref{betres} that:

\bl
\label{betresl}\mbox{}
\bi
\item[(i)] $\beta({\rm RRes}_{\ell}(T)) \subseteq {\rm RRes}_{\ell}(T)$.
\item[(ii)] $\beta({\rm LRes}_{\ell}(T)) \subseteq {\rm LRes}_{\ell}(T)$.
\ei
\el

The semaphore code $\RSC_{\ell}(T)$ consists of all minimal elements of $\rres_{\ell}(T)$ for the suffix order. Equivalently,
\beq
\label{prima2}
\RSC_{\ell}(T) = (\RSC(T) \cap \Omega^{\leq \ell}) \cup \Omega(\Omega^{\ell -1}\setminus \rres(T)).
\eeq
Dually, the left semaphore code $\LSC_{\ell}(T)$ consists of all minimal elements of $\lres_{\ell}(T)$ for the prefix order, or equivalently
\beq
\label{prima3}
\LSC_{\ell}(T) = (\LSC(T) \cap \Omega^{\leq \ell}) \cup (\Omega^{\ell -1}\setminus \lres(T))\Omega.
\eeq
Therefore $\RSC_{\ell}(T) \cup \LSC_{\ell}(T) \subseteq \Omega^{\leq \ell}$.

\subsection{A context-free example}

\subsubsection{Description}

We present now a very elementary example just to illustrate the notation and ideas. Further research will include much more complicated examples.

Let $A = \{ a,b \}$ and $L = \{ a^nb^n \mid n \geq 1\}$, one of the classical examples of a (deterministic) context-free language which is not rational. The language $L$ is accepted by the Turing machine $T$ depicted by
$$
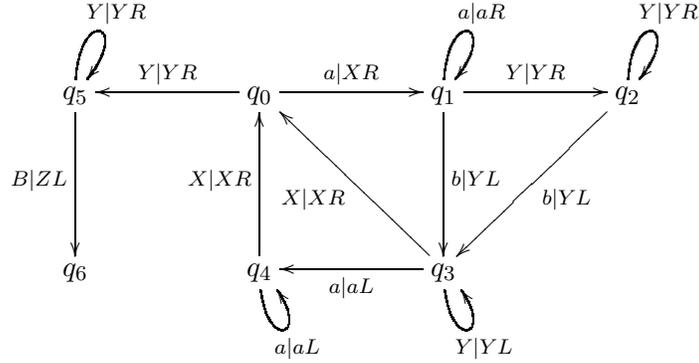
\xymatrix{
q_5 \ar@(u,ur)^{Y|{YR}} \ar[dd]_{B|{ZL}} && q_0 \ar[rr]^{a|{XR}} \ar[ll]_{Y|{YR}} && q_1 \ar[rr]^{Y|{YR}} \ar@(u,ur)^{a|{aR}} 
\ar[dd]^{b|{YL}} && q_2 \ar@(u,ur)^{Y|{YR}} \ar[ddll]^{b|{YL}} \\
&&&&&&\\
q_6 && q_4 \ar@(d,dr)_{a|{aL}} \ar[uu]^{X|{XR}} && q_3 \ar@(d,dr)_{Y|{YL}} \ar[ll]^{a|{aL}} \ar[uull]^{X|{XR}} && 
}$$
where $q_0$ is the initial state and $q_6$ the unique final state.

In state $q_0$ we can only read $a$ or $Y$. In the first case, $a$ is replaced by $X$ and we change to state $q_1$. 
Then we move right across other possible $a$'s until we reach the first $b$ and replace it by $Y$ to go to state $q_3$. 
If we have done this routine before, we may have to move across older $Y$'s -- taking us into state $q_2$. From state 
$q_3$, we intend to move left until we reach $X$, which means going through $Y$'s and then $a$'s (if there are some 
left -- state $q_4$). So we are back at state $q_0$ and we repeat the procedure. If we have replaced all the $a$'s, we 
are supposed to read $Y$ at state $q_0$, then move to the right end of the tape (state $q_5$) reading only $Y$'s. If 
we have succeeded on reaching the blank $B$, then we accept the input moving to the final state $q_6$.

It is easy to check that $T$ is legal-halting. Indeed, $B$ can be read at most once, and any long enough sequence of 
transitions must necessarily involve replacing $a$ by $X$ or $b$ by $Y$, which are both irreversible changes. 

We use the notation introduced in Section \ref{geco} for an arbitrary Turing machine.

\subsubsection{Resets}
\label{prima4}

We claim that 
\beq
\label{cfe1}
\begin{array}{lll}
\displaystyle\Omega^* \setminus \rres(T)&=&a^*\{ b,Y\}^* \cup (\bigcup_{i = 1,3,4} \, a^*a^{q_i}a^*\{ b,Y\}^*) \cup (\bigcup_{i = 1,2}\, a^*Y^*b^{q_i}\{ b,Y\}^*)\\
&&\\
&\cup&(\bigcup_{i = 1,2,3}\, a^*Y^*Y^{q_i}\{ b,Y\}^*).
\end{array}
\eeq

Let $m,n \geq 0$ and $u \in \{ b,Y\}^*$ containing $n$ $b$'s. Then 
$$\beta^{(\omega)}(a^{q_0}a^n\cdot a^mu,b,1) = Y \neq b = \beta^{(\omega)}(a^mu,b,1),$$
hence $a^mu \notin \rres(T)$.  

Assume now that $i \in \{ 1,3,4\}$, $m,n,k \geq 0$ and $u \in \{ b,Y\}^*$ contains $k$ $b$'s. Then 
$$\beta^{(\omega)}(Xa^{k+1}\cdot a^ma^{q_i}a^nu \cdot b,b,1) = Y \neq b = \beta^{(\omega)}(a^ma^{q_i}a^nu \cdot b,b,1),$$
hence $a^ma^{q_i}a^nu \notin \rres(T)$.

Assume next that $i \in \{ 1,2\}$, $m,n,k \geq 0$ and $u \in \{ b,Y\}^*$ contains $k$ $b$'s. Then 
$$\beta^{(\omega)}(Xa^{k+1}\cdot a^mY^nb^{q_i}u,b,1) = Y \neq b = \beta^{(\omega)}(a^mY^nb^{q_i}u,b,1),$$
hence $a^mY^nb^{q_i}u \notin \rres(T)$.

Finally, assume that $i \in \{ 1,2,3\}$, $m,n,k \geq 0$ and $u \in \{ b,Y\}^*$ contains $k$ $b$'s. Then 
$$\beta^{(\omega)}(Xa^{k+2}\cdot a^mY^nY^{q_i}u\cdot b,b,1) = Y \neq b = \beta^{(\omega)}(a^mY^nY^{q_i}u\cdot b,b,1),$$
hence $a^mY^nY^{q_i}u \notin \rres(T)$.

To prove the converse, we must prove that any other word is necessarily a right reset. We make extensive use from $\rres(T)$ being an ideal of $\Omega^*$.

Consider first 
$$r \in \{ B,X,Z \} \cup \{ B^q,X^q,Z^q \mid q \in Q \}.$$
Suppose that $t_1,t'_1,t_2,t_3 \in \Omega^*$, $P \in \Omega$ are such that both $t_1rt_2Pt_3, t'_1rt_2Pt_3 \in \leg(T)$. Let $n \geq 0$. If $\beta^{(n)}(t_1rt_2,P,t_3) \neq P$, then it is easy to see that $t_1 \in \Gamma^*$ and has no influence in the computation. Thus $\beta^{(n)}(t_1rt_2,P,t_3) = \beta^{(n)}(t'_1rt_2,P,t_3)$ and so $r \in \rres(T)$.

Thus we only need to discuss words $w \in \Omega^*$ such that $\tape(w) \in \{ a,b,Y\}^*$. We consider next the word $ba$. Consider $t_1bat_2Pt_3 \in \leg(T)$. If  $t_2Pt_3 \notin \Gamma^*$, then $t_1$ is irrelevant to the computation of $\beta^{(n)}(t_1bat_2,P,t_3)$. If $t_2Pt_3 \in \Gamma^*$, then $\beta^{(n)}(t_1bat_2,P,t_3) = P$ necessarily. It follows that $ba \in \rres(T)$. 

Similarly, $Ya, b^qa, ba^q, Y^qa, Ya^q \in \rres(T)$ for every $q \in Q$. Hence we have reduced the problem to words $w \in \Omega^*$ such that $\tape(w) \in a^*\{b,Y\}^*$.

Since the transition function is undefined for these pairs, we have
$$\{ a^{q_2}, a^{q_5}, a^{q_6}, b^{q_0}, b^{q_3}, b^{q_4}, b^{q_5}, b^{q_6}, Y^{q_4}, Y^{q_6} \} \subseteq \rres(T).$$
Also $a^{q_0} \in \rres(T)$ because $T$ moves to the right and will never get to the left of the new $X$. Similarly, $Y^{q_0}, Y^{q_5} \in \rres(T)$. To complete the proof of (\ref{cfe1}), it suffices to show that
\beq
\label{vic}
bY^*b^{q_i} \cup bY^*Y^{q_j} \subseteq \rres(T)
\eeq
for $i = 1,2$ and $j = 1,2,3$, because any word not containing such a factor has already been established to be or not to be a right reset.

Indeed, in neither case the head of $T$ can pass to the left of the first $b$, so (\ref{vic}) and therefore (\ref{cfe1}) hold as claimed.

Similarly, we compute the left resets, in fact we obtain
\beq
\label{sy}
\lres(T) = \rres(T).
\eeq

\subsubsection{Semaphore codes}
\label{prima5}

In view of (\ref{cfe1}), it is straightforward to check that
$$\begin{array}{lll}
\RSC(T)&=&(\Omega \setminus \{ a, a^{q_1}, a^{q_3}, a^{q_4}\})a^+\{ b,Y\}^*\\
&&\\ 
&\cup&
(\Omega \setminus \{ a, b, Y, a^{q_1}, a^{q_3}, a^{q_4}, b^{q_1},b^{q_2}, Y^{q_1}, Y^{q_2}, Y^{q_3} \})\{ b,Y\}^*\\
&&\\
&\cup& 
(\bigcup_{i = 1,3,4} \, (\Omega \setminus \{ a \})
a^*a^{q_i}a^*\{ b,Y\}^*) \cup 
(\bigcup_{i = 1,2}\, (\Omega \setminus \{ a \})a^+Y^*b^{q_i}\{ b,Y\}^*)\\
&&\\
&\cup&(\bigcup_{i = 1,2}\, (\Omega \setminus \{ a,Y \})Y^*b^{q_i}\{ b,Y\}^*)
\cup(\bigcup_{i = 1,2,3}\, (\Omega \setminus \{ a \})a^+Y^*Y^{q_i}\{ b,Y\}^*)\\
&&\\
&\cup&(\bigcup_{i = 1,2,3}\, (\Omega \setminus \{ a,Y \})Y^*Y^{q_i}\{ b,Y\}^*).
\end{array}$$

To compute the intersection $\RSC(T) \cap \leg(T)$, we replace the 5 last occurrences of $\Omega$ by $\Gamma$.

Similarly,
$$\begin{array}{lll}
\LSC(T)&=&
a^*(\Omega \setminus \{ a, b, Y, a^{q_1}, a^{q_3}, a^{q_4}, b^{q_1},b^{q_2}, Y^{q_1}, Y^{q_2}, Y^{q_3} \})\\
&&\\
&\cup&a^*Y^+(\Omega \setminus \{ b, Y, b^{q_1},b^{q_2}, Y^{q_1}, Y^{q_2}, Y^{q_3} \})
\cup a^*Y^*b\{ b,Y\}^*(\Omega \setminus \{ b, Y\})\\
&&\\
&\cup&(\bigcup_{i = 1,3,4} \, a^*a^{q_i}a^*(\Omega \setminus \{ a, b, Y\}))
\cup(\bigcup_{i = 1,3,4} \, a^*a^{q_i}a^*\{ b,Y\}^+(\Omega \setminus \{ b, Y\}))\\
&&\\
&\cup&
(\bigcup_{i = 1,2}\, a^*Y^*b^{q_i}\{ b,Y\}^*(\Omega \setminus \{ b, Y\}))
\cup
(\bigcup_{i = 1,2,3}\, a^*Y^*Y^{q_i}\{ b,Y\}^*(\Omega \setminus \{ b, Y\})).
\end{array}$$

To compute the intersection $\LSC(T) \cap \leg(T)$, we replace the 4 last occurrences of $\Omega$ by $\Gamma$.

\subsubsection{Semaphore codes modulo $\ell$}

In view of~\eqref{prima2} and~\eqref{prima3}, we can easily compute easily $\RSC_{\ell}(T)$ and $\LSC_{\ell}(T)$ 
making use of the computations performed in Sections~\ref{prima4} and~\ref{prima5}.

\section{Free pro-$\D$ semigroups}
\label{section.free profinite}

For general background on free pro-$\D$ semigroups, 
see~\cite[Sections 3.1 and 3.2]{RS:2009}.

Let $A$ be a finite nonempty alphabet. We denote by $A^{-\omega}$ the
set of all \defn{left infinite} words on $A$, that is, infinite sequences
of the form $\cdots a_3a_2a_1$ with $a_i \in A$. If $u \in A^+$, we
denote the left infinite word $\cdots uuu$ by $u^{-\omega}$.

The free semigroup $A^+$ acts on the right of $A^{-\omega}$ by
concatenation: given $x = \cdots a_3a_2a_1 \in A^{-\omega}$ and $u =
a'_1\cdots a'_n \in A^+$, we define
$$xu = \cdots a_3a_2a_1a'_1\cdots a'_n \in A^{-\omega}.$$
Given $x \in A^+ \cup A^{-\omega}$ and $y \in A^{-\omega}$, we define
also $xy = y$. Together with concatenation on $A^+$, this defines a
semigroup structure for $A^+ \cup A^{-\omega}$. 

The \defn{suffix (ultra)metric} on $A^+ \cup A^{-\omega}$ is defined as
follows. Given $x,y \in A^+ \cup A^{-\omega}$, let
$\lcs(x,y)$ be the longest common suffix of $x$ and $y$, and define
$$d(x,y) = \left\{
\begin{array}{ll}
2^{-|\lcs(x,y)|}&\mbox{ if }x \neq y,\\
0&\mbox{ otherwise.}
\end{array}
\right.$$
Given $x_0 \in A^+ \cup A^{-\omega}$ and $\delta > 0$, we write
$$B_{\delta}(x_0) = \{ x \in A^+ \cup A^{-\omega} \mid d(x,x_0) < \delta \}$$
for the open ball of radius $\delta$ around $x_0$.

If $S \in \D$ is endowed with the discrete topology and $\p\colon A \to S$
is a mapping, then there exists a unique continuous homomorphism 
$\Phi\colon A^+ \cup A^{-\omega} \to S$ such that the diagram
$$\xymatrix{
A_{}  \ar[r]^{\p} \ar@{^{(}->}[d] & S\\
A^+ \cup A^{-\omega} \ar@{-->}[ur]_{\Phi} & \  
}$$
commutes. This characterizes $(A^+ \cup A^{-\omega},d)$ as the
\defn{free pro-$\D$ semigroup} on $A$. We shall denote it by
$\oo{\Omega}_A(\D)$. It is well known that $\oo{\Omega}_A(\D)$ is
a complete and compact topological semigroup.
 
We remark that for a general pseudovariety $\V$, the metric considered
for free pro-$\V$ semigroups is the profinite metric, but in the
particular case of $\D$ we can use this alternative metric that
equates to the normal form.

\section{$(-\omega)$-reset graphs}
\label{section.regr2}

We consider now $A$-graphs with possibly infinite vertex sets. 
For general concepts in automata theory, the reader is referred to~\cite{Berstel:1979}.

A \defn{left infinite path} in an $A$-graph $\Gamma = (Q,E)$ is an infinite
sequence of the form
$$\cdots \mapright{a_3} q_3 \mapright{a_2} q_2 \mapright{a_1} q_1$$
such that $(q_{i+1},a_i,q_i) \in E$ and $q_i\in Q$ for every $i \geq 1$. Its label is
the left infinite word $\cdots a_3a_2a_1 \in A^{-\omega}$. We write 
$$\cdots \mapright{x} q$$
to denote a left infinite path with label $x$ ending at $q$. 

An $A$-graph $\Gamma = (Q,E)$ is:
\bi
\item
\defn{deterministic} if
$(p,a,q), (p,a,q') \in E$ implies $q = q'$;
\item
\defn{complete} if for all $p \in Q$ and $a \in A$ there exists some edge $(p,a,q) \in E$;
\item
\defn{strongly connected} if, for all $p,q \in Q$, there exists a path
$p \mapright{u} q$ in $\Gamma$ for some $u \in A^*$;
\item
\defn{$(-\omega)$-deterministic} if
$$\cdots \mapright{x} q,\hspace{.3cm} \cdots \mapright{x} q'
\mbox{ paths in } \Gamma \Rw q = q'$$
holds for all $q,q' \in Q$ and $x \in A^{-\omega}$;
\item
\defn{$(-\omega)$-complete} if every $x \in A^{-\omega}$ labels some left
infinite path in $\Gamma$;
\item
\defn{$(-\omega)$-trim} if every $q \in Q$ occurs in some left
infinite path in $\Gamma$;
\item
a \defn{$(-\omega)$-reset graph} if it is $(-\omega)$-deterministic,
$(-\omega)$-complete and $(-\omega)$-trim.
\ei
We denote by $\rg(A)$ the class of all $(-\omega)$-reset $A$-graphs.

If $\Gamma = (Q,E) \in \rg(A)$, then $Q$ induces a partition
$$A^{-\omega} = \bigcup_{q \in Q} A^{-\omega}_q,$$
where $A^{-\omega}_q$ denotes the set of all $x \in A^{-\omega}$
labelling some path $\cdots \mapright{x} q$ in $\Gamma$. Moreover,
$A^{-\omega}_q \neq\emptyset$ for every $q \in Q$. 

\bp
\label{deco}
Let $\Gamma \in {\rm RG}(A)$. Then $\Gamma$ is deterministic and
complete.
\ep

\proof
Write $\Gamma = (Q,E)$ and suppose that $(p,a,q),(p,a,q') \in
E$. Since $\Gamma$ is $(-\omega)$-trim, there exists some left infinite
path $\cdots \mapright{x} p$ for some $x \in A^{-\omega}$. Hence there
exist left infinite
paths $\cdots \mapright{xa} q$ and $\cdots \mapright{xa} q'$, and
since $\Gamma$ is $(-\omega)$-deterministic, we get $q = q'$. Therefore
$\Gamma$ is deterministic.

Let $p \in Q$ and $a \in A$. Since $\Gamma$ is $(-\omega)$-trim, there
exists some left infinite 
path $\cdots \mapright{x} p$ for some $x \in A^{-\omega}$. Now $xa \in
A^{-\omega}$ and $\Gamma$ being $(-\omega)$-complete implies that there
exists some path $\cdots \mapright{xa} q$ in $\Gamma$, which we may
factor as
$$\cdots \mapright{x} q' \mapright{a} q.$$
Since $\Gamma$ is $(-\omega)$-deterministic, we get $q' = p$, hence
$(p,a,q) \in E$ and $\Gamma$ is complete.
\qed

We recall now the preorder $\leq$ introduced in~\cite[Section 3]{RSS:2016}.
Given $A$-graphs $\Gamma,\Gamma'$, we write $\Gamma \leq \Gamma'$ if
there exists a morphism $\Gamma \to \Gamma'$.

\bl
\label{newpaor}
Let $A$ be a finite nonempty alphabet and let $\Gamma,\Gamma' \in
{\rm RG}(A)$ with  $\Gamma \leq \Gamma' \leq \Gamma$. Then $\Gamma 
\cong \Gamma'$.   
\el

\proof
Let $\p:\Gamma \to \Gamma'$ and $\psi:\Gamma' \to \Gamma$ be morphisms. 
Write $\Gamma = (Q,E)$ and $\Gamma' = (Q',E')$. It is easy to see that
$$A^{-\omega}_q \subseteq A^{-\omega}_{q\p}, \quad A^{-\omega}_{q'}
\subseteq A^{-\omega}_{q'\psi}$$
for all $q \in Q$ and $q' \in Q'$. Hence $A^{-\omega}_q \subseteq
A^{-\omega}_{q\p\psi}$. Since $\Gamma$ is $(-\omega)$-trim, we have
$A^{-\omega}_q \neq \emptyset$. Since $\Gamma$ is
$(-\omega)$-deterministic, we get $q = q\p\psi$. Similarly, $q' =
q'\psi\p$, hence $\p$ and $\psi$ are mutually inverse
bijections and therefore mutually inverse $A$-graph isomorphisms.
\qed

Let $[\Gamma]$ denote the isomorphism class of $\Gamma$.
Similarly to~\cite[Section 3]{RSS:2016},
$$[\Gamma] \leq [\Gamma'] \hspace{.5cm}\mbox{if } \Gamma \leq \Gamma'$$
defines a preorder on $\rg(A)/\cong$. Moreover, Lemma \ref{newpaor} yields:

\bc
\label{newporg}
Let $A$ be a finite nonempty alphabet. Then $\leq$
is a partial order on ${\rm RG}(A) /\cong$.
\ec

\section{Right congruences on $A^{-\omega}$}
\label{section.rc omega}

Since $xy = y$ for all $x \in \oo{\Omega}_A(\D)$ and $y \in A^{-\omega}$,
it follows that $A^{-\omega}$ is the minimum ideal of $\oo{\Omega}_A(\D)$. 
Following the notation introduced in~\cite[Section 2.2]{RSS:2016}, we denote by $\rc(A^{-\omega})$ 
the lattice of right congruences on $A^{-\omega}$ (with respect to the right action
of $\oo{\Omega}_A(\D)$). 

We say that $\rho \in \rc(A^{-\omega})$ is \defn{closed} if $\rho$ is a
closed subset of $A^{-\omega} \times A^{-\omega}$ for the product metric
$$d'((x,y),(x',y')) = \max\{ d(x,x'),d(y,y') \},$$
where $d$ denotes the suffix metric on $A^{-\omega}$. 
Given $x_0,y_0 \in A^+ \cup A^{-\omega}$ and $\delta > 0$, we write
$$B_{\delta}((x_0,y_0)) = \{ (x,y) \in (A^+ \cup A^{-\omega})^2 \mid d'((x,y),(x_0,y_0)) < \delta \}.$$
By~\cite[Exercise 3.1.7]{RS:2009}, this implies that $x\rho$ is a closed
subset of $A^{-\omega}$ for every $x \in A^{-\omega}$. The next
example shows that the converse fails.

\be
\label{cnc}
Let $A = \{ a,b\}$ and let 
\beq
\label{cnc1}
w = \ldots a^4ba^3ba^2bab.
\eeq
For all $x,y \in A^{-\omega}$, let
$$x\rho y \quad \text{ if } \quad \left\{
\begin{array}{l}
x = wu,\; y = wv \mbox{ with } |u| = |v|\\
\hspace{.5cm}\mbox{or}\\
x = y.
\end{array}
\right.$$
Then $\rho \in {\rm RC}(A^{-\omega})$ and $x\rho$ is closed for every
$x \in A^{-\omega}$, but $\rho$ is not closed.
\ee

Indeed, it is easy to see that, given $x \in A^{-\omega}$, there is at
most one word $u \in A^*$ such that $x = wu$. We call this a
\defn{$w$-factorization} of $x$. Hence $\rho$ is transitive
and it follows immediately that $\rho \in \rc(A^{-\omega})$. The
uniqueness of the $w$-factorization implies also that $x\rho$ is
finite (hence closed) 
for every $x \in A^{-\omega}$. However,
$$\lim_{n\to \infty} (wa^n,wb^n) = (a^{-\omega},b^{-\omega}) \notin \rho.$$
Since $(wa^n,wb^n) \in \rho$ for every $n \geq 1$, then $\rho$ is not
closed.

\medskip

We denote by $\crc(A^{-\omega})$ (respectively $\orc(A^{-\omega})$)
the set of all closed (respectively open) right congruences on
$A^{-\omega}$. 

We consider $\crc(A^{-\omega})$ (partially) ordered by
inclusion. Similarly to~\cite[Section 5]{RSS:2016}, we can relate
$\crc(A^{-\omega})$ with $\rg(A)$. 

Given $\rho \in \rc(A^{-\omega})$, the Cayley graph $\cay(\rho)$ is
the $A$-graph $\cay(\rho) = (A^{-\omega}/\rho,E)$ defined by
$$E = \{ (u\rho,a,(ua)\rho) \mid u \in A^{-\omega}, \; a \in A \}.$$

\bl
\label{lip}
Let $\rho \in {\rm RC}(A^{-\omega})$.
\bi
\item[(i)] For every $x \in A^{-\omega}$,
there exists a left infinite path $\cdots \mapright{x} x\rho$ in
$\cay(\rho)$.
\item[(ii)] ${\rm Cay}(\rho)$ is $(-\omega)$-complete and $(-\omega)$-trim.
\ei 
\el

\proof
(i) Write $x = \cdots a_3a_2a_1$ with $a_i \in A$. 
For every $n
\geq 1$, write $x_n = \cdots a_{n+2}a_{n+1}a_n$. Then 
$$\cdots \mapright{a_3} x_3\rho \mapright{a_2} x_2\rho \mapright{a_1}
x_1\rho = x\rho$$
is a left infinite path in $\cay(\rho)$ labeled by $x$.

(ii) By part (i).
\qed

\bl
\label{propcay}
Let $\rho \in {\rm CRC}(A^{-\omega})$. Then:
\bi
\item[(i)] If $\cdots \mapright{x} q$ is a left infinite path in
  ${\rm Cay}(\rho)$, then $q = x\rho$.
\item[(ii)] ${\rm Cay}(\rho) \in {\rm RG}(A)$.
\ei
\el

\proof
(i) Assume that $q = y\rho$ with $y \in
A^{-\omega}$. Write $x = \cdots a_3a_2a_1$ with $a_i \in A$. 
For every $n
\geq 1$, let $u_n =
a_n \cdots a_1$.
Then there exists some path 
$y_{n}\rho \mapright{u_n} y\rho$ in $\cay(\rho)$ for some $y_n \in
A^{-\omega}$. Hence $y_{n}u_n \in
y\rho$. Since 
$$x = \lim_{n\to \infty} u_n = \lim_{n\to \infty} y_{n}u_n$$ and $\rho$ 
closed implies $y\rho$ closed, we get $x \in y\rho$, hence $x\rho =
y\rho = q$.

(ii) By part (i), $\cay(\rho)$ is $(-\omega)$-deterministic. By Lemma
\ref{lip}(ii), 
$\cay(\rho)$ is both $(-\omega)$-complete and $(-\omega)$-trim, therefore 
$\cay(\rho) \in \rg(A)$.
\qed

We discuss next open right congruences, relating them in particular
with the right congruences on $A^k$. 
Given $x \in A^{-\omega}$, let $x\xi_k$ denote the suffix of length $k$ of $x$.
For $\sigma \in \rc(A^k)$, let
$\wh{\sigma}$ be the relation on $A^{-\omega}$ defined by
$$x \wh{\sigma} y \hspace{.3cm} \mbox{ if } (x\xi_k) \sigma
(y\xi_k).$$
It is immediate that $\wh{\sigma} \in \rc(A^{-\omega})$. 

On the other hand, given $\rho \in \rc(A^{-\omega})$ and $k \geq 1$, we
define a relation $\rho^{(k)}$ on $A^k$ by
$$u \rho^{(k)} v \hspace{.3cm} \mbox{ if }(A^{-\omega}u \times
A^{-\omega}v) \cap \rho \neq \emptyset.$$
We denote by $\rho^{[k]}$ the transitive closure of $\rho^{(k)}$.

The next example shows that $\rho^{(k)}$ needs not to be
transitive, even in the closed case. 

\be
\label{notr}
Let $A = \{ a,b\}$ and let $w$ be given by (\ref{cnc1}). For all $x,y
\in A^{-\omega}$, let 
$$x\rho y \quad \text{ if } \quad \left\{
\begin{array}{l}
\{ x, y\} = \{ wa^2u, wbau\} \mbox{ for some }u \in A^*\\
\hspace{.5cm}\mbox{or}\\
\{ x, y\} = \{ wb^2av, wb^3v\} \mbox{ for some }v \in A^*\\
\hspace{.5cm}\mbox{or}\\
x = y.
\end{array}
\right.$$
Then $\rho \in {\rm CRC}(A^{-\omega})$ but $\rho^{(2)}$ is not transitive.
\ee

Indeed, by the uniqueness of the $w$-factorization remarked in Example
\ref{cnc}, $\rho$ turns out to be transitive and therefore a right
congruence. 

We sketch the proof that $\rho$ is closed. Let $(x,y) \in (A^{-\omega} \times
A^{-\omega}) \setminus \rho$. Then $x \neq y$. Write $u =
\lcs(x,y)$. We consider several cases:

\smallskip
\noindent
\underline{Case I}: $\{ x,y \} = \{ zb^2au,z'b^3u \}$.

\smallskip
\noindent
Then either $z \neq w$ or $z' \neq w$. We may assume that $z \neq w$. Let
$k \geq 1$ be such that $w \notin B_{2^{-k}}(z)$. It is easy to see
that $B_{2^{-k-3-|u|}}((x,y)) \cap \rho = \emptyset$. 

\smallskip
\noindent
\underline{Case II}: $\{ x,y \} = \{ za^2u,z'bau \}$.

\smallskip
\noindent
Then either $z \neq w$ or $z' \neq w$. We may assume that $z \neq w$. Let
$k \geq 1$ be such that $w \notin B_{2^{-k}}(z)$. It is easy to see
that $B_{2^{-k-2-|u|}}((x,y)) \cap \rho = \emptyset$. 

\smallskip
\noindent
\underline{Case III}: all the remaining cases.

\smallskip
\noindent
It is easy to see that $B_{2^{-3-|u|}}((x,y)) \cap \rho = \emptyset$. 

Therefore $\rho$ is closed.

Now $(wa^2,wba) \in \rho$ yields $(a^2,ba) \in
\rho^{(2)}$, and $(wb^2a,wb^3) \in \rho$ yields $(ba,b^2) \in
\rho^{(2)}$, However, $(a^2,b^2) \notin \rho^{(2)}$, hence
$\rho^{(2)}$ is not transitive. 

\medskip

The following lemma compiles some elementary properties of
$\rho^{(k)}$ and $\rho^{[k]}$. The proof is left to the reader.

\bl
\label{elpr}
Let $A$ be a finite nonempty alphabet, $\rho \in {\rm
  RC}(A^{-\omega})$ and $k \geq 1$. Then:
\bi
\item[(i)] $\rho^{(k)} \in {\rm RC}(A^k)$ if and only if
$\rho^{(k)}$ is transitive;
\item[(ii)] $\rho^{[k]} \in {\rm RC}(A^k)$;
\item[(iii)] $\rho \subseteq \ds\bigcap_{n \geq 1} \wh{\rho^{[n]}}$.
\ei
\el

We discuss next some alternative characterizations for open right
congruences. We recall the definition of \defn{$k$-reset graph} from \cite[Section 3]{RSS:2016}. 

We say that $u \in A^*$ is a \defn{reset word} for a deterministic
and complete $A$-graph $\Gamma = (Q,E)$ if
$|Qu| = 1$. This is equivalent to say that all paths labeled by $u$
end at the same vertex. Let $\res(\Gamma)$ denote the set of all reset
words for $\Gamma$. 

We say that $\Gamma$ is a \defn{$k$-reset graph} if $A^k \subseteq \res(\Gamma)$. We denote by $\rg_k(A)$ 
the class of all strongly connected deterministic complete $k$-reset $A$-graphs.

\bp
\label{open}
Let $A$ be a finite nonempty alphabet and 
$\rho \in {\rm RC}(A^{-\omega})$. Then the following conditions are equivalent:
\bi
\item[(i)] $\rho$ is open;
\item[(ii)] $x\rho$ is an open subset of $A^{-\omega}$ for every
  $x \in A^{-\omega}$; 
\item[(iii)] $\rho = \wh{\sigma}$ for some $\sigma \in {\rm RC}(A^k)$
  and $k \geq 1$;
\item[(iv)] there exists some $k \geq 1$ such that $\rho^{(k)}$ is
  transitive and $\rho = \wh{\rho^{(k)}}$;
\item[(v)] ${\rm Cay}(\rho) \in {\rm RG}_k(A)$ for some $k \geq 1$;
\item[(vi)] $\rho$ is closed and has finite index. 
\ei
\ep

\proof
(i) $\Rw$ (ii). Let $x \in A^{-\omega}$. Since
$(x,x) \in \rho$, there exists some $\delta > 0$ such that
$B_{\delta}((x,x)) \subseteq \rho$. Since $B_{\delta}((x,x)) =
B_{\delta}(x) \times B_{\delta}(x)$, we get
$B_{\delta}(x) \subseteq x\rho$ and so $x\rho$ is open.

(ii) $\Rw$ (vi). Let $x,y \in A^{-\omega}$ be
such that $(x,y) \notin \rho$. Since $x\rho$ and $y\rho$ are open,
there exists some $\delta 
> 0$ such that $B_{\delta}(x)
\subseteq x\rho$ and $B_{\delta}(y) \subseteq y\rho$. If $x' \in
x\rho$ and $y' \in y\rho$, then $(x,y) \notin \rho$ yields $(x',y')
\notin \rho$. Hence
$$(B_{\delta}(x) \times
B_{\delta}(y)) \cap \rho = \emptyset,$$
and so $B_{\delta}((x,y)) \cap \rho = \emptyset$. Thus the
complement of $\rho$ is open and so $\rho$ is closed.

On the other hand, $\{ x\rho \mid x \in A^{-\omega} \}$ is an open
cover of $A^{-\omega}$ and so admits a finite subcover since
$A^{-\omega}$ is compact. Therefore $\rho$ has finite index.

(vi) $\Rw$ (v). By Lemma \ref{propcay}(ii), we have $\cay(\rho) \in
\rg(A)$. Hence $\cay(\rho)$ is deterministic and complete by Proposition
\ref{deco}. 

Let $x,y \in A^{-\omega}$. Since $\cay(\rho)$ is $(-\omega)$-trim,
there exists a left infinite path $\cdots \mapright{z} y\rho$ in
$\cay(\rho)$. Since $\rho$ has finite index, we may factor this path
as
$$\cdots \mapright{z'} w\rho \mapright{u} w\rho \mapright{v} y\rho$$
with $u \neq \varepsilon$. On the other hand, since $\cay(\rho)$ is
complete and $\rho$ has finite index, there exist $m \geq 0$ and $p
\geq 1$ such that there exists a path
$$x\rho \mapright{u^m} x'\rho \mapright{u^p} x'\rho$$
in $\cay(\rho)$. It follows that there exist two paths
$$\cdots \mapright{u^{-\omega}} w\rho, \quad 
\cdots \mapright{u^{-\omega}} x'\rho$$
and so $w\rho = x'\rho$ since $\cay(\rho)$ is
$(-\omega)$-deterministic. Thus there exists a path
$$x\rho \mapright{u^m} w\rho \mapright{v} y\rho$$
and so $\cay(\rho)$ is strongly connected.

Suppose now that $\cay(\rho) \notin \rg_k(A)$ for every $k \geq
1$. Let $P$ denote the set of pairs of distinct vertices in
$\cay(\rho)$. Then
$$\forall k\geq 1\; \exists u_k \in A^k\; \exists (p,q) \in P\;
\exists \mbox{ paths } \cdots \mapright{u_k} p,\; \cdots
\mapright{u_k} q \mbox{ in }\cay(\rho).$$
Since $P$ is finite, one of the pairs $(p,q)$ must repeat infinitely
often. Hence there exists some $(p,q) \in P$ such that
$$\forall k\geq 1\; \exists u_k \in A^{\geq k}\;
\exists \mbox{ paths } \cdots \mapright{u_k} p,\; \cdots
\mapright{u_k} q \mbox{ in }\cay(\rho).$$
Since $\oo{\Omega}_A(\D)$ is compact, we may replace $(u_k)_k$ by some
convergent subsequence. Let $x = \lim_{k \to \infty} u_k$. Since
$(|u_k|)_k$ is unbounded, we have $x \in A^{-\omega}$.

Write $p = x_p\rho$ with $x_p \in A^{-\omega}$. 
Since $\cay(\rho)$ is $(-\omega)$-trim, there exists some left
infinite path $\cdots \vlongmapright{y_ku_k} x_p\rho$ for some $y_k \in
A^{-\omega}$. By Lemma \ref{lip}(i), there exists a path
$\cdots \vlongmapright{y_ku_k} (y_ku_k)\rho$ in $\cay(\rho)$. Since
$\cay(\rho)$ is $(-\omega)$-deterministic, we get $(y_ku_k)\rho =
x_p\rho$, hence $y_ku_k \in
x_p\rho$. Since $\rho$ closed implies $x_p\rho$ closed and 
$$x = \lim_{k \to \infty} u_k = \lim_{k \to \infty} y_ku_k,$$
we get $x \in x_p\rho$. By Lemma \ref{lip}(i), there exists a path
$\cdots \mapright{x} x\rho = x_p\rho = p$ in $\cay(\rho)$. Similarly,
there exists some path $\cdots \mapright{x} q$. Since $p \neq q$, this
contradicts $\cay(\rho)$ being $(-\omega)$-deterministic. Therefore
$\cay(\rho) \in \rg_k(A)$ for some $k \geq 1$.

(v) $\Rw$ (iv). 
Assume that $\cay(\rho) \in \rg_k(A)$ for some $k \geq
1$. We show that
\beq
\label{open1}
x\xi_k = y\xi_k \Rw x\rho y
\eeq
holds for all $x,y \in A^{-\omega}$. Indeed, by Lemma \ref{lip}(i), there
exists left infinite paths
$$\cdots \mapright{x} x\rho, \quad \cdots \mapright{y}
y\rho$$
in $\cay(\rho)$. Since $x\xi_k = y\xi_k \in A^k \subseteq \res(\cay(\rho))$, we get
$x\rho = y\rho$ and so (\ref{open1}) holds.

Suppose now that $u,v,w \in A^k$ are such that $u \rho^{(k)} v
\rho^{(k)} w$. Then there exist some $x,y,y',z \in A^{-\omega}$ such
that $(xu)\rho(yv)$ and $(y'v)\rho(zw)$. Then $(yv)\xi_k = v =
(y'v)\xi_k$ and (\ref{open1}) yields $(yv)\rho(y'v)$. Thus
$(xu)\rho(zw)$ by transitivity and so  
$u \rho^{(k)} w$. Therefore $\rho^{(k)}$ is transitive.

Now it follows from Lemma \ref{elpr} that $\wh{\rho^{(k)}}$ is well
defined and $\rho \subseteq \wh{\rho^{(k)}}$.

Conversely, let $(x,y) \in \wh{\rho^{(k)}}$. Then $(x\xi_k,y\xi_k) \in
\rho^{(k)}$ and so there exist $x',y' \in A^{-\omega}$ such that
$(x'(x\xi_k),y'(y\xi_k)) \in \rho$. Since $(x'(x\xi_k))\xi_k =
x\xi_k$, it follows from (\ref{open1}) that $(x'(x\xi_k))\rho
x$. Similarly, $(y'(y\xi_k))\rho y$ and we get $x\rho y$ by
transitivity. Therefore $\wh{\rho^{(k)}} \subseteq \rho$ as required.

(iv) $\Rw$ (iii). In view of Lemma \ref{elpr}(i).

(iii) $\Rw$ (i). Let $(x,y) \in \rho = \wh{\sigma}$ and let $(x',y')
\in B_{2^{-k}}((x,y))$. Then $x'\xi_k = x\xi_k$ and $y'\xi_k = y\xi_k$.
Hence
$$x\rho y \Rw (x\xi_k)\sigma(y\xi_k) \Rw (x'\xi_k)\sigma(y'\xi_k) \Rw
x'\rho y'$$
and so $B_{2^{-k}}((x,y)) \subseteq \rho$. Therefore $\rho$ is open.
\qed

The following example shows that closed is required in condition (vi).

\be
\label{cir}
Let $A = \{ a,b\}$ and let $\rho$ be the relation on $A^{-\omega}$
defined by $x\rho y$ if $b$ occurs in both $x,y$ or in none of
them. Then $\rho$ is a right congruence of index 2 on $A^{-\omega}$
but it is not closed.
\ee

Indeed, it is immediate that $\rho$ is a right congruence of index
2. Since $a^{-\omega} = \lim_{n\to \infty} b^{-\omega}a^n$, $\rho$ is
not closed.

\medskip

We say that $\rho \in \rc(A^{-\omega})$ is \defn{profinite} if $\rho$ is
an intersection of open right congruences. Since open right
congruences are closed by Proposition \ref{open}, it follows that
every profinite right congruence, being the intersection of closed
sets, is itself closed.
We denote by $\prc(A^{-\omega})$
the set of all profinite right congruences on $A^{-\omega}$.

Given a graph $\Gamma = (Q,E)$ and $k \geq 1$, we define a relation
$\mu^{(k)}_{\Gamma}$ on $Q$ by
$$p \mu^{(k)}_{\Gamma} q\hspace{.5cm}\mbox{if there exist paths
}\cdots \mapright{u} p,\; \cdots \mapright{u} q \mbox{ in
  $\Gamma$ for some }u \in A^k.$$ 
Let $\mu^{[k]}_{\Gamma}$ denote the reflexive and transitive closure
of $\mu^{(k)}_{\Gamma}$. Then $\mu^{[k]}_{\Gamma}$ is an equivalence
relation on $Q$.

\bp
\label{prof}
Let $A$ be a finite nonempty alphabet and $\rho \in {\rm
  RC}(A^{-\omega})$. Then the following conditions are equivalent:
\bi
\item[(i)] $\rho$ is profinite;
\item[(ii)] $\rho$ is an intersection of countably many open congruences;
\item[(iii)] $\rho = \ds\bigcap_{k \geq 1} \wh{\rho^{[k]}}$;
\item[(iv)] $\ds\bigcap_{k \geq 1} \mu^{[k]}_{{\rm Cay}(\rho)} = id$. 
\ei
\ep

\proof
(i) $\Rw$ (iii). Assume that $\rho = \cap_{i \in I} \tau_i$ with
$\tau_i \in \orc(A^{-\omega})$ for every $i \in I$.

We have $\rho \subseteq \cap_{k
  \geq 1} \wh{\rho^{[k]}}$ by Lemma \ref{elpr}(iii). To prove the opposite
inclusion, we show that
\beq
\label{prof1}
\forall i \in I\; \exists k \geq 1 \; \wh{\rho^{[k]}} \subseteq \tau_i.
\eeq

Indeed, it follows from Proposition \ref{open} that there exist some
$k \geq 1$ and $\sigma_i \in \rc(A^k)$ such that $\tau_i = 
\wh{\sigma_i}$. We claim that
\beq
\label{prof2}
\tau_i^{(k)} \subseteq \sigma_i.
\eeq
Assume that $(u,v) \in \tau_i^{(k)}$. Then there exist $x,y \in
A^{-\omega}$ such that $(xu,yv) \in \tau_i = \wh{\sigma_i}$. Hence
$$(u,v) = ((xu)\xi_k,(yv)\xi_k) \in \sigma_i$$
and (\ref{prof2}) holds.

Since $\rho \subseteq \tau_i$ implies $\rho^{(k)} \subseteq
\tau_i^{(k)}$, it follows that $\rho^{(k)} \subseteq \sigma_i$ and
so $\rho^{[k]} \subseteq \sigma_i$ since $\sigma_i$ is
transitive. Thus
$$\wh{\rho^{[k]}} \subseteq \wh{\sigma_i} = \tau_i$$
and (\ref{prof1}) holds. 

Therefore
$$\ds\bigcap_{k
  \geq 1} \wh{\rho^{[k]}} \subseteq \cap_{i \in I} \tau_i = \rho$$
as required.

(iii) $\Rw$ (ii). By Lemma \ref{elpr}(ii), $\rho^{[k]} \in \rc(A^k)$ for
every $k \geq 1$, hence $\wh{\rho^{[k]}}$ is open by Proposition
\ref{open} and we are done.

(ii) $\Rw$ (i). Trivial.

(iii) $\Rw$ (iv). Write $\mu^{(k)} = \mu^{(k)}_{{\rm Cay}(\rho)}$ and $\mu^{[k]}
= \mu^{[k]}_{{\rm Cay}(\rho)}$.  By Lemma \ref{elpr}(ii), $\rho^{[k]}
\in \rc(A^k)$ for every $k \geq 1$, hence $\wh{\rho^{[k]}}$ is open
(and therefore closed) by Proposition \ref{open}. Therefore $\rho$ is
closed and so $\cay(\rho) \in \rg(A)$ by Lemma \ref{propcay}(ii).

Let $x,y \in A^{-\omega}$ be such that $x\rho \neq
y\rho$. Suppose that $(x\rho,y\rho) \in \mu^{[k]}$. Then there exist
$z_0,\ldots,z_n \in A^{-\omega}$ such that $z_0 = x$, $z_n = y$ and
$(z_{i-1}\rho,z_i\rho) \in \mu^{(k)}$ for $i = 1,\ldots,n$. For $i =
1,\ldots,n$, there exist paths
$$z'_{i-1}\rho \mapright{u_i} z_{i-1}\rho,\quad z''_{i}\rho \mapright{u_i} z_i\rho$$
in $\cay(\rho)$ for some $u_i \in A^k$ and $z'_{i-1},z''_{i} \in A^{-\omega}$.

Hence $z_{i-1}\rho = (z'_{i-1}u_i)\rho$ and $z_i\rho =
(z''_{i}u_i)\rho$, yielding 
$$(z_{i-1}\xi_k)\,\rho^{(k)}\, u_i\, \rho^{(k)}\,
(z_{i}\xi_k)$$
and so $(z_{i-1}\xi_k)\rho^{[k]} (z_{i}\xi_k)$. Now
$(x\xi_k)\rho^{[k]} (y\xi_k)$ follows by transitivity, hence $(x,y)
\in \wh{\rho^{[k]}}$. Since $x\rho \neq y\rho$ implies $(x,y)
\notin \wh{\rho^{[m]}}$ for some $m \geq 1$ by condition (iii), it follows that 
$(x\rho,y\rho) \notin \mu^{[m]}$ and so (iv) holds.

(iv) $\Rw$ (iii). By Lemma \ref{elpr}(iii), we have $\rho \subseteq
\cap_{k \geq 1} \wh{\rho^{[k]}}$. Conversely, let $(x,y) \in
\cap_{k \geq 1} \wh{\rho^{[k]}}$. For each $k$, we have
$(x\xi_k,y\xi_k) \in \rho^{[k]}$, hence there exist $u_0, \ldots,u_n
\in A^k$ such that $u_0 = x\xi_k$, $u_n = y\xi_k$ and $(u_{i-1},u_i)
\in \rho^{(k)}$ for $i = 1,\ldots,n$. For $i = 1,\ldots,n$, there exist
$z_{i-1},z'_i \in A^{-\omega}$ such that $(z_{i-1}u_{i-1},z'_iu_i) \in
\rho$. Write also $x = z'_0u_0$ and $y = z_nu_n$.

By Lemma \ref{lip}(i), there exist paths
$$\cdots \longmapright{z'_iu_i} (z'_iu_i)\rho, \quad 
\cdots \longmapright{z_iu_i} (z_iu_i)\rho$$
in $\cay(\rho)$ for $i = 0,\ldots,n$, hence
$$((z_{i-1}u_{i-1})\rho,(z_iu_i)\rho) = ((z'_iu_i)\rho,(z_iu_i)\rho) \in
\mu^{(k)}.$$
Thus
$$((z_0u_0)\rho,(z_{n}u_{n})\rho) \in
\mu^{[k]}.$$
Since $((z'_0u_0)\rho,(z_{0}u_{0})\rho) \in
\mu^{(k)}$, we get 
$$(x\rho,y\rho) = ((z'_0u_0)\rho,(z_{n}u_{n})\rho) \in
\mu^{[k]}.$$
Since $k$
is arbitrary, it follows from condition (iv) that $x\rho = y\rho$,
hence $\cap_{k \geq 1} \wh{\rho^{[k]}} \subseteq \rho$ as required.
\qed

Every open right congruence on $A^{-\omega}$ is trivially profinite
and we remarked before that every profinite right congruence is
necessarily closed. Hence
$$\orc(A^{-\omega}) \subseteq \prc(A^{-\omega}) \subseteq
\crc(A^{-\omega}).$$
We show next that these inclusions are strict if $|A| > 1$.

For every $k \geq 1$, let $\rho_k$ be the relation on $A^{-\omega}$
defined by
$$x\rho_k y \hspace{.5cm}\mbox{if }x\xi_k = y\xi_k.$$
It is easy to check that $\rho_k \in \orc(A^{-\omega})$ for every $k
\geq 1$. Since $\cap_{k \geq 1} \rho_k = id$, it follows that the
identity congruence is profinite, while it is clearly not open.

To construct a closed non profinite right congruence is much
harder. We do it through the following example. 

\be
\label{cnp}
Let $A = \{ a,b\}$. Given $u,v \in A^k$, write $u < v$ if $u = u'aw$
and $v = v'bw$ for some $w \in A^*$. Let $u_1^{(k)} < \cdots <
u_{2^k}^{(k)}$ be the elements of $A^k$, totally ordered by $<$. Let
$p_1 < p_2 < \cdots$ be 
the prime natural numbers. For every $n \in \N$, let
$$w_{n} = \ldots a^{p_3^n}ba^{p_2^n}ba^{p_1^n}b.$$
Let $\rho \in {\rm RC}(A^{-\omega})$ be generated by the relation
$$R = \{ (w_{p_k^i}u_{i}^{(k)},w_{p_k^i}u_{i+1}^{(k)}) \mid k \geq 1,\; 1 \leq i
< 2^k \} \cup \{ (b^{-\omega}a,a^{-\omega}b) \}.$$
Then $\rho$ is closed but not profinite.
\ee

We start by showing that
\beq
\label{cnp1}
w_{p_k^i}A^* \cap w_{p_{k'}^{i'}}A^* \neq \emptyset \; \mbox{ implies
} \; (k = k' \mbox{ and } i = i')
\eeq
for all $k,k',i,i' \geq 1$. Indeed, suppose that $w_{p_k^i}u =
w_{p_{k'}^{i'}}v$ for some $u,v \in A^*$. By definition of $w_n$, $w_{p_k^i}u$
has only finitely many factors of the form $ba^{2m}b$, and the
leftmost must be $ba^{2^{p_k^i}}b$. Since $w_{p_k^i}u =
w_{p_{k'}^{i'}}v$, we get $ba^{2^{p_k^i}}b = ba^{2^{p_{k'}^{i'}}}b$
and so $p_k^i = p_{k'}^{i'}$. Therefore $k = k'$ and $i = i'$, and
(\ref{cnp1}) holds.

Write
$$R' = \{ (xu,yu) \mid (x,y) \in R \cup R\inv,\; u \in A^* \}.$$
Let $x \in A^{-\omega}$. We show that 
\beq
\label{cnp2}
\mbox{there exists at most one $y \in
A^{-\omega}$ such that $(x,y) \in R'$}. 
\eeq
This is obvious if $x \in
b^{-\omega}aA^* \cup a^{-\omega}bA^*$, hence we may assume that $x \in
w_{p_k^i}A^*$ for some $k \geq 1$ and $1 \leq i < 2^k$. In
view of (\ref{cnp1}), we must have
$$\{ x,y \} = \{ w_{p_k^i}u_{i}^{(k)}v, w_{p_k^i}u_{i+1}^{(k)}v \}$$
for some $v \in A^*$, and $k,i$ are uniquely determined. Since
$w_{p_k^i} \notin w_{p_k^i}A^+$, also $u_{i}^{(k)}$, $u_{i+1}^{(k)}$
and $v$ are uniquely determined. Thus (\ref{cnp2}) holds.

Suppose that $x\, R' \, y \, R'\,z$ with $x \neq y \neq z$. Since $R'$
is symmetric, (\ref{cnp2}) yields $x' = z'$. It follows that $R' \cup
\; id$ is an equivalence relation, indeed the smallest right congruence
containing $R$. It follows that
$$R' \cup \; \mathrm{id} = \rho.$$
Moreover, each $\rho$-class contains at most two
elements. 

We prove now that $\rho$ is closed. Let $(x,y) \in (A^{-\omega} \times
A^{-\omega}) \setminus \rho$. Then $x \neq y$, hence we may assume
without loss of generality that $x = x'av$ and $y = y'bv$ with $v \in
A^*$. Let
$$m = \max\{ i \geq 0 \mid b^i <_s x',\; a^i <_s y'\}.$$
Note that the above set is bounded, otherwise $x' = b^{-\omega}$ and
$y' = a^{-\omega}$, yielding 
$$(x,y) = (b^{-\omega}av,a^{-\omega}bv) \in \rho,$$
a contradiction. Write $x' = x''b^m$ and $y' = y''a^m$.

For $j = 0,\ldots,|v|$, write $v = v_jv'_j$ with $|v_j| = j$. Then
$a^mbv_j$ is the successor of $b^mav_j$ in the ordering of
$A^{m+1+j}$, hence we may write
$$b^mav_j = u_{i_j}^{(m+1+j)},\quad a^mbv_j = u_{i_j+1}^{(m+1+j)}$$
for some $1 \leq i_j < 2^{m+1+j}$. It follows that
\beq
\label{cnp5}
x = x''u_{i_j}^{(m+1+j)}v'_j,\quad y = y''u_{i_j+1}^{(m+1+j)}v'_j
\eeq
Let
$$m_j = \min\{ |\lcs(x'', w_{p_{m+1+j}^{i_j}})|, |\lcs(y'',
w_{p_{m+1+j}^{i_j}})| \}.$$
Note that $m_j$ is a well-defined natural number, otherwise $x'' =
w_{p_{m+1+j}^{i_j}} = y''$ and 
$$(x,y) =
(w_{p_{m+1+j}^{i_j}}u_{i_j}^{(m+1+j)}v'_j,w_{p_{m+1+j}^{i_j}}u_{i_j+1}^{(m+1+j)}v'_j)
\in \rho,$$
a contradiction. 

Let
$$p = \max\{ m_0, \ldots, m_{|v|} \} +m+1+|v|.$$
We show that
\beq
\label{cnp3}
B_{2^{-p}}((x,y)) \cap \rho = \emptyset.
\eeq

Suppose that $(z_1,z_2) \in B_{2^{-p}}((x,y)) \cap \rho$. Since $p >
1+|v|$, we have $av <_s z_1$ and $bv <_s z_2$. By
maximality of $m$, and since $p > m+1+|v|$, we have either
$ab^mav <_s z_1$ or $ba^mbv <_s z_2$. Hence we must have
\beq
\label{cnp4}
z_1 = w_{p_k^i}u_{i}^{(k)}v', \quad z_2 = w_{p_k^i}u_{i+1}^{(k)}v'
\eeq
for some $v'$, $k \geq 1$ and $1 \leq i < 2^k$. Clearly, $|v'| \leq |v|$,
hence we must have $v' = v'_j$ for $j = |v|-|v'|$. 

We have $x = x''b^mav_jv'_j$, hence $b^mav_jv'_j <_s z_1$. Similarly,
$y = y''a^mbv_jv'_j$ yields $a^mbv_jv'_j <_s z_2$.
Suppose that $k < m+1+j$. Since $|\lcs(x,z_1)| > m+1+|v|$, it follows
from (\ref{cnp4}) that $w_{p_k^i}$ ends with a $b$. Similarly,
$|\lcs(y,z_2)| > m+1+|v|$ implies that $w_{p_k^i}$ ends with an $a$, a
contradiction. Hence $k \geq m+1+j$.

Suppose now that $k > m+1+j$. By maximality of $m$, we must have one of
the following cases:
\bi
\item
$ab^mav_j \leq_s u_{i}^{(k)}$ and $a^mbv_j \leq_s u_{i+1}^{(k)}$;
\item
$b^mav_j \leq_s u_{i}^{(k)}$ and $ba^mbv_j \leq_s u_{i+1}^{(k)}$.
\ei
Any of these cases contradicts $u_{i+1}^{(k)}$ being the successor of
$u_{i}^{(k)}$ for the ordering of $A^k$, hence $k = m+1+j$ and we may
write
$$z_1 = w_{p_{m+1+j}^i}u_{i}^{(m+1+j)}v'_j, \quad z_2 =
w_{p_{m+1+j}^i}u_{i+1}^{(m+1+j)}v'_j.$$
Since $d(z_1,x) < 2^{-m_j-m-1-|v'|} = 2^{-m_j-m-1-j-|v'_j|}$, it
follows from (\ref{cnp5}) that $i = i_j$ and 
$$|\lcs(x'',w_{p_{m+1+j}^{i_j}})| > m_j.$$
Similarly,  
$$|\lcs(y'',w_{p_{m+1+j}^{i_j}})| > m_j,$$
contradicting the definition of $m_j$.

Thus (\ref{cnp3}) holds and so $(A^{-\omega} \times
A^{-\omega}) \setminus \rho$ is open. Therefore $\rho$ is closed.

Let $k \geq 1$. Since $(w_{p_k^i}u_{i}^{(k)},w_{p_k^i}u_{i+1}^{(k)})
\in \rho$, we have $(u_{i}^{(k)},u_{i+1}^{(k)}) \in \rho^{(k)}$ for
every $1 \leq i < 2^k$. Since $u_{1}^{(k)} = a^k$ and $u_{2^k}^{(k)} =
b^k$, it follows that
$(a^k,b^k) \in \rho^{[k]}$
and so $(a^{-\omega},b^{-\omega}) \in \wh{\rho^{[k]}}$. Since $k$ is
arbitrary, we get
$$(a^{-\omega},b^{-\omega}) \in \ds\bigcap_{k\geq 1} \wh{\rho^{[k]}}.$$
However, 
$$(a^{-\omega},b^{-\omega}) \notin R' \cup \; \mathrm{id} = \rho,$$
hence $\rho \neq \cap_{k\geq 1} \wh{\rho^{[k]}}$ and so $\rho$ is not
profinite by Proposition \ref{prof}.

\section{Special right congruences on $A^{-\omega}$}
\label{section.src omega}

To avoid trivial cases, we assume throughout this section that $A$ is
a finite alphabet containing at least two elements. 

Given $P \subseteq A^*$, we define a relation $\tau_P$ on $A^{-\omega}$ by:
$$x \tau_P y \; \mbox{ if $x = y$ or $x,y \in A^{-\omega}u$ for some }u \in P.$$
 
\bl
\label{newirco}
Let $P \subseteq A^*$. Then $\tau_P$ is an equivalence relation on $A^{-\omega}$.
\el

\proof
It is immediate that $\tau_P$ is reflexive and symmetric. For
transitivity, we may assume that $x,y,z \in A^{-\omega}$ are distinct
and $x \, \tau_P \, y  \, \tau_P \, z$. Then there exist $u,v \in P$
such that $u <_s x,y$ and $v <_s y,z$. Since $u$ and $v$ are
both suffixes of $y$, one of them is a suffix of the other. Hence
either $u <_s x,z$ or $v <_s x,z$. Therefore $\tau_P$ is
transitive.
\qed

If we consider left ideals, being a right congruence turns out to be a
special case: 

\bp
\label{newnsrc}
Let $L \unlhd_{\ell} A^*$. Then the following conditions are equivalent:
\bi
\item[(i)] $\tau_L \in {\rm RC}(A^{-\omega})$;
\item[(ii)] $\tau_L \in {\rm PRC}(A^{-\omega})$;
\item[(iii)] $L \unlhd A^*$;
\item[(iv)] $(L\beta_{\ell})A \subseteq A^*(L\beta_{\ell})$.
\item[(v)] $L\beta_{\ell}$ is a semaphore code.
\ei
\ep

\proof
(i) $\Rw$ (iv). We may assume that $|A| > 1$.
Let $u \in L$ and $a \in A$. Take $b \in
A \setminus \{ a \}$. Then $(a^{-\omega}u,b^{-\omega}u) \in \tau_L$
and by (i) we get $(a^{-\omega}ua,b^{-\omega}ua) \in \tau_L$. It
follows that $ua$ has some suffix in $L$, hence $LA \subseteq A^*L =
L$ and so 
$$(L\beta_{\ell})A \subseteq LA \subseteq L = A^*(L\beta_{\ell}).$$

(iv) $\Rw$ (iii). We have
$$LA = A^*(L\beta_{\ell})A \subseteq A^*(L\beta_{\ell}) = L.$$
It follows that $LA^*\subseteq L$. Since $L \unlhd_{\ell} A^*$, we get $L
\unlhd A^*$. 

(iii) $\Rw$ (ii). By Lemma \ref{newirco}, $\tau_L$ is an equivalence
relation.
Let $x,y \in A^{-\omega}$ be such that $x\tau_L y$. We
may assume that there exists some $u \in L$ such that $u <_s x,y$.
Since $L \unlhd A^*$, we have $ua \in L$ and $ua <_s xa,ya$ yields
$(xa,ya) \in \tau_L$. Therefore $\tau_L \in
\rc(A^{-\omega})$.

Let $(x,y) \in (A^{-\omega} \times A^{-\omega}) \setminus \tau_L$. Then
$x \neq y$. Let $u = \lcs(x,y)$ and let $m = |u|+1$. We claim that
\beq
\label{newnsrc1} 
(x\xi_m)\tau_L^{[m]} = \{ x\xi_m \}.
\eeq
Indeed, suppose that $(x\xi_m,v) \in \tau_L^{(m)}$ and $v \neq
x\xi_m$. Then there exist $z,z' \in A^{-\omega}$ such that
$(z(x\xi_m),z'v) \in \tau_L$. Since $v \neq
x\xi_m$, then $z(x\xi_m)$ and $z'v$ must have a common suffix $w \in L$,
and $|w| < m$. But then $w \leq_s u$, yielding $u \in L$ and
$x\tau_Ly$, a contradiction. Thus $(x\xi_m,v) \in \tau_L^{(m)}$ implies $v =
x\xi_m$, and so (\ref{newnsrc1}) holds.

Suppose that $(x,y) \in \wh{\tau_L^{[m]}}$. Then $(x\xi_m,y\xi_m)
\in \tau_L^{[m]}$, hence $x\xi_m = y\xi_m$ by (\ref{newnsrc1}),
contradicting $m > |u|$. Thus $(x,y) \notin \wh{\tau_L^{[m]}}$ and so
$\cap_{k \geq 1} \wh{\tau_L^{[m]}}  \subseteq \tau_L$. Hence $\tau_L =
\cap_{k \geq 1} \wh{\tau_L^{[m]}}$ by Lemma \ref{elpr}(iii), and so
$\tau_L$ is profinite by Proposition \ref{prof}.

(ii) $\Rw$ (i). Trivial.

(iv) $\iff$ (v). By Lemma~\cite[Lemma 4.1]{RSS:2016}, since $L\beta_{\ell}$ is always
a suffix code.
\qed

We say that $\rho \in \rc(A^{-\omega})$ is a \defn{special right
  congruence} on $A^{-\omega}$ if 
$\rho = \tau_I$ for some $I \unlhd A^*$. 
In view of Proposition~\ref{newnsrc}, this is equivalent to say that $\rho = \tau_S$ for 
some semaphore code $S$ on $A$. We denote by  
$\src(A^{-\omega})$ the set of all special right congruences on $A^{-\omega}$.

The next result characterizes the open special right
congruences. Recall that a suffix code $S \subset A^*$ is said to be
\defn{maximal} 
if $S \cup \{ u \}$ fails to be a suffix code for every $u \in A^*
\setminus S$.

\bp
\label{osrc}
Let $I \unlhd A^*$. Then the following conditions are equivalent:
\bi
\item[(i)] $\tau_I \in {\rm ORC}(A^{-\omega})$;
\item[(ii)] $I\beta_{\ell}$ is a finite maximal suffix code;
\item[(iii)] $A^* \setminus I$ is finite.
\ei
\ep

\proof
(i) $\Rw$ (ii). Let $u,v \in I\beta_{\ell}$ be distinct. Then
$$((A^{-\omega}u) \times (A^{-\omega}v)) \cap \tau_I = \emptyset.$$
Since $\tau_I$ has finite index by Proposition \ref{open}, it follows
that the suffix code $I\beta_{\ell}$ is finite.

Suppose now that $I\beta_{\ell} \cup \{ u \}$ is a suffix code for
some $u \in A^* \setminus (I\beta_{\ell})$. It is easy to
see that no two elements of $A^{-\omega}u$ are $\tau_I$ equivalent, a
contradiction since $\tau_I$ has finite index. Therefore
$I\beta_{\ell}$ is a maximal suffix code.

(ii) $\Rw$ (iii). Let 
$m$ denote the maximum length of the words in $I\beta_{\ell}$. Suppose
that $v \in A^* \setminus I$ has length $> m$. It is straightforward
to check that $I\beta_{\ell} \cup \{ v \}$ is a suffix code,
contradicting the maximality of $I\beta_{\ell}$. Thus $A^* \setminus I
\subseteq A^{\leq m}$ and is therefore finite.

(iii) $\Rw$ (i). We have $\tau_I \in \rc(A^{-\omega})$ by Proposition
\ref{newnsrc}. 

Let $m \geq 1$ be such that $A^* \setminus I
\subseteq A^{\leq m}$. Then
$$x\xi_{m+1} = y\xi_{m+1} \; \Rw \; x\tau_Iy$$
holds for all $x,y \in A^{-\omega}$ and so $\tau_I$ has finite
index. Since $\tau_I$ is profinite (and therefore closed) by
Proposition \ref{newnsrc}, it follows from Proposition \ref{open} that
$\tau_I$ is open.
\qed

The proof of~\cite[Lemma 7.4]{RSS:2016} can be adapted to show that inclusion
among left ideals determines inclusion for the equivalence relations
$\tau_L$: 

\bl
\label{newordi}
Let $|A| > 1$ and $L,L' \unlhd_{\ell} A^*$. Then 
$$\tau_L \subseteq \tau_{L'} \iff L \subseteq L'.$$
\el

Note that Lemma \ref{newordi} does not hold for $|A| = 1$, since
$|A^{-\omega}| = 1$. 

Similarly, we adapt~\cite[Proposition 7.6]{RSS:2016}:

\bp
\label{newisolat}
Let $|A| > 1$. Then:
\bi
\item[(i)] $\tau_{I\cap J} = \tau_I \cap \tau_J$ and $\tau_{I\cup J} =
  \tau_I \cup \tau_J$ for all $I,J \unlhd A^*$;
\item[(ii)] ${\rm SRC}(A^{-\omega})$ is a full sublattice of ${\rm
    RC}(A^{-\omega})$; 
\item[(iii)] the mapping
$$\begin{array}{rcl}
\I(A)&\to&{\rm SRC}(A^{-\omega})\\
I&\mapsto&\tau_I
\end{array}$$
is a lattice isomorphism.
\ei
\ep

Given $\rho \in \rc(A^{-\omega})$ and $C \in A^{-\omega}/\rho$, we say that
$C$ is \defn{nonsingular} if $|C| > 1$. If $C$ is nonsingular, we
denote by $\lcs(C)$ the 
longest common suffix of all words in $C$. We define
\bi
\item $\Lambda_{\rho} = \{ \lcs(C) \mid C \in A^{-\omega}/\rho \mbox{
    is nonsingular} \},$
\item $\Lambda'_{\rho} = \{ \lcs(x,y) \mid (x,y) \in \rho, \; x \neq y \}$.
\ei

\bl
\label{newpropsl}
Let $\rho \in {\rm RC}(A^{-\omega})$. Then:
\bi
\item[(i)] $A^*\Lambda_{\rho} = A^*\Lambda'_{\rho}$;
\item[(ii)] $\Lambda'_{\rho} \unlhd_r A^*$; 
\item[(iii)] $A^*\Lambda'_{\rho} \unlhd A^*$. 
\ei
\el

\proof
(i) Let $C \in A^{-\omega}/\rho$ be nonsingular and let $w = \lcs(C)$. By
maximality of $w$ there exist $a,b 
\in A$ distinct and $x,y \in A^{-\omega}$ such that $xaw,ybw \in C$. Thus $w =
\lcs(xaw,ybw)$ and so 
\beq
\label{newpropsl1}
\Lambda_{\rho} \subseteq
\Lambda'_{\rho}.
\eeq
Therefore $A^*\Lambda_{\rho} \subseteq
A^*\Lambda'_{\rho}$.

Conversely, let $(x,y) \in \rho$ with $x \neq y$. Then $x\rho$ is
nonsingular and $\lcs(x\rho)$
is a suffix of 
$\lcs(x,y)$, hence $\Lambda'_{\rho} \subseteq
A^*\Lambda_{\rho}$ and so $A^*\Lambda_{\rho} =
A^*\Lambda'_{\rho}$.

(ii) Let $u \in \Lambda'_{\rho}$ and $a \in A$. Then $u = \lcs(x,y)$
for some $(x,y) 
\in \rho$ with $x \neq y$. Then $(xa,ya) 
\in \rho$. Since $\lcs(xa,ya) = ua$, we get $ua \in
\Lambda'_{\rho}$. Therefore $\Lambda'_{\rho} \unlhd_r A^*$.

(iii) Clearly, $A^*\Lambda'_{\rho} \unlhd_{\ell} A^*$. Now we use part
(ii).
\qed

Given $\rho \in \rc(A^{-\omega})$, we write
$$\res(\rho) = \res(\cay(\rho)).$$
We refer to the elements of $\res(\rho)$ as the resets of $\rho$.

\bl
\label{newprops}
Let $\rho \in {\rm RC}(A^{-\omega})$. Then:
\bi
\item[(i)] $\res(\rho) \unlhd A^*$;
\item[(ii)] if $\rho$ is closed, then 
$$\res(\rho) = \{ w \in A^* \mid
  (xw,yw) \in \rho \mbox{ for all }x,y \in A^{-\omega} \}.$$
\ei
\el

\proof
(i) Immediate.

(ii) Let $w \in \res(\rho)$ and $x,y \in A^{-\omega}$. By Lemma
\ref{lip}(i), there exist paths 
$$\cdots \mapright{xw} (xw)\rho, \quad \cdots \mapright{yw} (yw)\rho$$
in $\cay(\rho)$. Now $w \in \res(\rho)$ yields $(xw)\rho = (yw)\rho$.

Now let $w \in A^* \setminus \res(\rho)$. Then there exist paths $p
\mapright{w}q$ and $p'
\mapright{w}q'$ in $\cay(\rho)$ with $q \neq q'$. Since $\cay(\rho)$
is $(-\omega)$-trim by Lemma
\ref{lip}(ii), there exist left infinite paths 
$$\cdots \mapright{x} p, \quad \cdots \mapright{y} p'$$
in $\cay(\rho)$, hence paths
$$\cdots \mapright{xw} q, \quad \cdots \mapright{yw} q'.$$
Since $\rho$ is closed, it follows from Lemma \ref{propcay}(i) that
$(xw)\rho = q \neq q' = (yw)\rho$ and we are done.
\qed

Adapting the proof of~\cite[Proposition 7.9]{RSS:2016}, we obtain:

\bp
\label{newinclu}
Let $|A| > 1$, $\rho \in {\rm RC}(A^{-\omega})$ and $I \unlhd A^*$. Then:
\bi
\item[(i)] $\rho \subseteq \tau_I \iff \Lambda_{\rho} \subseteq I \iff
\Lambda'_{\rho} \subseteq I$;
\item[(ii)] if $\rho$ is closed, then $\tau_I \subseteq \rho \iff I
  \subseteq \res(\rho)$. 
\ei
\ep

Given $R \subseteq A^{-\omega}
 \times A^{-\omega}$, we denote by $R^{\sharp}$ the right
congruence on $A^{-\omega}$ generated by $R$, i.e. the intersection of all right
congruences on $A^{-\omega}$ containing $R$.

Given $\rho \in \rc(A^{-\omega})$, we denote by
$\ns(\rho)$ the set of all nonsingular $\rho$-classes.

We can now prove several equivalent characterizations of special right
congruences. 

\bp
\label{newsrc}
Let $|A| > 1$ and $\rho \in {\rm RC}(A^{-\omega})$. Then the following
conditions are equivalent:
\bi
\item[(i)] $\rho \in {\rm SRC}(A^{-\omega})$;
\item[(ii)] ${\rm lcs}:{\rm NS}(\rho) \to A^*$ is injective,
  $\Lambda_{\rho}$ is a suffix code and
\beq
\label{extr}
\forall x \in A^{-\omega}\; \forall w \in \Lambda_{\rho}\; (xw)\rho
\in {\rm NS}(\rho);
\eeq
\item[(iii)] $\rho = \tau_{A^*\Lambda_{\rho}}$;
\item[(iv)] $\rho = \tau_{A^*\Lambda'_{\rho}}$;
\item[(v)] $\rho = \tau_L^{\sharp}$ for some $L \unlhd_{\ell} A^*$.
\ei
\ep

\proof
(i) $\Rw$ (ii).
By a straightforward adaptation of the proof of (i) $\Rw$ (ii) 
in~\cite[Proposition 7.10]{RSS:2016}, we check that ${\rm lcs}:{\rm NS}(\rho) \to
A^*$ is injective and $\Lambda_{\rho}$ is a suffix code.

Now let $x \in A^{-\omega}$ and $w \in \Lambda_{\rho}$. Then $w =
\lcs(y\rho)$ for some $y\rho \in \ns(\rho)$. By (\ref{newpropsl1}), we
may write $w = \lcs(y',y'')$ for some $y',y'' \in y\rho$
distinct. Since $\rho = \tau_I$, it follows that $w \in I$, hence
$(xw,y),(xw,y') \in \tau_I = \rho$. Since $y' \neq y''$, it follows
that either $xw \neq y'$ or $xw \neq y''$, so in any case $(xw)\rho
\in \ns(\rho)$ as required.

(ii) $\Rw$ (iii). Write $I = A^*\Lambda_{\rho}$. 
If $(x,y) \in \rho$ and $x \neq y$, then $\lcs(x\rho) \in
\Lambda_{\rho} \subseteq I$ 
is a suffix of both $x$ and $y$, hence $(x,y) \in \tau_I$. Thus $\rho
\subseteq \tau_I$.

Conversely, let $(x,y) \in \tau_I$. We may assume that $x \neq y$,
hence there exists some $w \in \Lambda_{\rho}$
such that $w <_s x,y$. Hence (\ref{extr}) yields $x\rho,y\rho \in
\ns(\rho)$. 

Suppose that $\lcs(x\rho) \neq w$. Then
$\lcs(x\rho) <_s w$ or $w <_s \lcs(x\rho)$, contradicting
$\Lambda_{\rho}$ being a suffix code. Hence $\lcs(x\rho) =
w$. Similarly, $\lcs(y\rho) = w$. Since $\lcs:\ns(\rho) \to
A^*$ is injective, we get $x\rho = y\rho$. Thus $\rho = \tau_I$.

(iii) $\iff$ (iv). By Lemma \ref{newpropsl}(i).

(iii) $\Rw$ (v). Write $L = A^*\Lambda_{\rho}$. By (iii), we have $\tau_L^{\sharp}
= \rho^{\sharp} = \rho$. Since $L \unlhd A^*$ by Lemma \ref{newpropsl},
(iv) holds.

(v) $\Rw$ (i). Let $I = LA^* \unlhd A^*$. Since $L \subseteq I$, it
follows from Lemma \ref{newordi} that $\tau_L \subseteq \tau_I$, hence 
$$\rho = \tau_L^{\sharp} \subseteq \tau_I^{\sharp} = \tau_I$$
by Proposition \ref{newnsrc}. 

Conversely, let $(x,y) \in \tau_I$. We may assume that $x \neq
y$. Then there exist factorizations $x 
= x'w$ and $y = y'w$ with $w \in I$. Write $w = zw'$ with $z \in
L$. Then $(x'z,y'z) \in \tau_L$ and so
$$(x,y) = (x'w,z'w) = (x'zw',y'zw') \in \tau_L^{\sharp} = \rho.$$  
Thus $\tau_I \subseteq \rho$ as required.
\qed

\bp
\label{newcsrc}
Let $|A| > 1$ and $\rho \in {\rm CRC}(A^{-\omega})$. Then the following
conditions are equivalent:
\bi
\item[(i)] $\rho \in {\rm SRC}(A^{-\omega})$;
\item[(ii)] ${\rm lcs}:{\rm NS}(\rho) \to A^*$ is injective,
  $\Lambda_{\rho}$ is a suffix code and
$$\forall x \in A^{-\omega}\; \forall w \in \Lambda_{\rho}\; (xw)\rho
\in {\rm NS}(\rho);$$
\item[(iii)] $\rho = \tau_{A^*\Lambda_{\rho}}$;
\item[(iv)] $\rho = \tau_{A^*\Lambda'_{\rho}}$;
\item[(v)] $\rho = \tau_L^{\sharp}$ for some $L \unlhd_{\ell} A^*$;
\item[(vi)] $\rho = \tau_{\res(\rho)}$;
\item[(vii)] $\Lambda_{\rho} \subseteq \res(\rho)$;
\item[(viii)] $\Lambda'_{\rho} \subseteq \res(\rho)$;
\item[(ix)] whenever
\beq
\label{newsrc3}
p \mapright{aw} q,\quad p' \mapright{bw} q,\quad p'' \mapright{w}
r
\eeq
are paths in ${\rm Cay}(\rho)$ with $a,b \in A$ distinct, then $q = r$.
\ei
\ep

\proof
(i) $\iff$ (ii) $\iff$ (iii) $\iff$ (iv) $\iff$ (v). By Proposition
\ref{newsrc}. 

(i) $\Rw$ (vi). If $\rho = \tau_I$ for some $I \unlhd A^*$, then 
$I \subseteq \res(\rho)$ by Proposition \ref{newinclu}(ii). Since
$\res(\rho) \unlhd A^*$ by Lemma \ref{newprops}(i), then Proposition
\ref{newinclu}(ii) also yields 
$$\tau_{\res(\rho)} \subseteq \rho = \tau_I,$$
hence $\res(\rho) \subseteq I$ by Lemma \ref{newordi}. Therefore $I =
\res(\rho)$.

(vi) $\Rw$ (vii) $\iff$ (viii). By Lemma \ref{newprops}(i), $\res(\rho)
\unlhd A^*$. Now we apply Proposition \ref{newinclu}(i).

(viii) $\Rw$ (i). We have $A^*\Lambda'_{\rho},\res(\rho) \unlhd A^*$
by Lemmas \ref{newpropsl}(iii) and \ref{newprops}(i).
It follows from Proposition \ref{newinclu} that
$$\tau_{\res(\rho)} \subseteq \rho \subseteq
\tau_{A^*\Lambda'_{\rho}}.$$
Since $\Lambda'_{\rho} \subseteq \res(\rho)$ yields
$A^*\Lambda'_{\rho} \subseteq \res(\rho)$ and therefore
$\tau_{A^*\Lambda'_{\rho}} \subseteq \tau_{\res(\rho)}$ by Lemma
\ref{newordi}, we get $\rho = \tau_{\res(\rho)} \in \src(A^{-\omega})$.

(viii) $\Rw$ (ix). Consider the paths in (\ref{newsrc3}). 
By Lemma \ref{lip}(ii), there exist left infinite paths
$$\cdots \mapright{x} p, \quad \cdots \mapright{x'} p'$$
in $\cay(\rho)$, hence $(xaw,x'bw) \in \rho$ by Lemma
\ref{propcay}(i) and so
$$w = \lcs(xaw,x'bw) \in \Lambda'_{\rho} \subseteq \res(\rho).$$
Thus $q = r$ as required.

(ix) $\Rw$ (viii). Let $w \in \Lambda'_{\rho}$. Then $w = \lcs(x,y)$
for some $(x,y) \in \rho$ such that $x \neq y$. 
We may write $x = x'aw$ and $y
= y'bw$ with $a,b \in A$ distinct. By Lemma \ref{lip}(i), there exist in 
$\cay(\rho)$ paths of the form 
$$\cdots \mapright{x'} p \mapright{aw} x\rho,\quad \cdots \mapright{y'} p'
\mapright{bw} y\rho = x\rho.$$ 
Now (ix) implies that $w \in \res(\rho)$.
\qed

We can now prove that not all open right congruences are special, even for
$|A| = 2$:

\be
\label{newnotsp}
Let $A = \{ a,b\}$ and let $\sigma$ be the equivalence relation on $A^3$
defined by the following partition:
$$\{ a^3,aba,ba^2\} \cup \{ bab,a^2b \} \cup \{ ab^2 \} \cup \{ b^2a
\} \cup \{ b^3 \}.$$
Then $\wh{\sigma} \in {\rm ORC}(A^3) \setminus {\rm SRC}(A^3)$.
\ee

Indeed, it is routine to check that $\sigma \in \rc(A^3)$, hence
$\rho = \wh{\sigma} \in {\rm ORC}(A^{-\omega})$ by Proposition \ref{open}. 
Since
$\lcs(a^{-\omega}\rho) = a$ and $\lcs((b^{-\omega}a)\rho) = b^2a$,
then $\Lambda_{\rho}$ 
is not a suffix code and so $\rho \notin \src(A^{-\omega})$ by
Proposition~\ref{newsrc}. 

\medskip

Let $\rho \in \rc(A^{-\omega})$ and let 
$$\underline{\rho} = \vee\{ \tau \in \src(A^{-\omega}) \mid \tau \subseteq
\rho \},$$
$$\oo{\rho} = \wedge\{ \tau \in \src(A^{-\omega}) \mid \tau \supseteq
\rho \}.$$
By Proposition \ref{newisolat}(ii), we have $\underline{\rho}, \oo{\rho}
\in \src(A^{-\omega})$.

\bp
\label{newmima}
Let $|A| > 1$ and $\rho \in {\rm CRC}(A^{-\omega})$. Then:
\bi
\item[(i)] $\underline{\rho} = \tau_{\res(\rho)}$;
\item[(ii)] $\oo{\rho} = \tau_{A^*\Lambda_{\rho}} = \tau_{A^*\Lambda'_{\rho}}$.
\ei
\ep

\proof
(i) By Lemma \ref{newprops}(i), we have $\res(\rho) \unlhd A^*$. Now the
claim follows from Proposition \ref{newinclu}(ii).

(ii) Similarly, we have $A^*\Lambda_{\rho} =
A^*\Lambda'_{\rho} \unlhd A^*$ by Lemma \ref{newpropsl}(iii), and the
claim follows from Proposition \ref{newinclu}(i).
\qed

The straightforward adaptation of~\cite[Example 7.14]{RSS:2016} shows that the pair
$(\underline{\rho},\oo{\rho})$ does not univocally determine $\rho \in
\rc(A^{-\omega})$, even in the open case:

\be
\label{newcer}
Let $A = \{ a,b\}$ and let $\sigma,\sigma'$ be the equivalence
relations on $A^3$ 
defined by the following partitions:
$$\{ a^3,aba,ba^2\} \cup \{ bab,a^2b \} \cup \{ ab^2 \} \cup \{ b^2a
\} \cup \{ b^3 \},$$
$$\{ a^3,b^2a,ba^2\} \cup \{ bab,a^2b \} \cup \{ ab^2 \} \cup \{ aba
\} \cup \{ b^3 \}.$$
Let $\rho = \wh{\sigma}$ and $\rho' = \wh{\sigma'}$. Then $\rho,\rho'
\in {\rm ORC}(A^{-\omega})$, $\underline{\rho} = 
\underline{\rho'}$ and $\oo{\rho} = \oo{\rho'}$.
\ee

This same example shows also that $\oo{\rho}$ does not necessarily
equal or cover $\underline{\rho}$ in $\src(A^{-\omega})$. Indeed, in
  this case we have 
$$\res(\rho) = A^*A^3 \cup \{ a^2,ab\} \subset I \subset A^+ \setminus
\{ b,b^2 \} = A^*\Lambda_{\rho}$$ 
for $I = A^*A^3 \cup \{
a^2,ab,ba\} \unlhd A^*$. By Lemma \ref{newordi}, we get
$$\underline{\rho} \subset \tau_I \subset \oo{\rho}.$$

\section{Conclusion and future work}
\label{section.future}

We enter now into random walks on infinite semigroups. The most sophisticated approach is described in~\cite{HM:2011}. 
We use profinite limits (see~\cite{RS:2009}) as an alternative approach, as developed in 
Sections~\ref{section.free profinite}-\ref{section.src omega}.

Indeed, if $I_1,I_2, \ldots$ is a sequence of ideals in $A^*$ with $I = \cap_{k \geq 1} I_{k}$, let $J_k\beta_{\ell}$ the 
semaphore code determined by the ideal $J_{k} = I_1 \cap \ldots \cap I_{k}$. Whenever $k \geq m$, we may define 
a mapping $\p_{km} \colon J_k\beta_{\ell} \to J_m\beta_{\ell}$ by setting $u\p_{km}$ to be the unique suffix of $u$ in 
$J_m\beta_{\ell}$. It is routine to check that:
\bi
\item
$\p_{km}$ is onto;
\item
$\p_{km}$ preserves the action of $A^*$ on the right;
\item
$(\p_{km})$ constitutes a projective system of surjective morphisms with respect to this action;
\item
$I\beta_{\ell}$ is the projective limit of this system.
\ei

In view of~\eqref{prima1}, each Turing machine $T$ provides an instance of this setting when $I_k = \rres_k(T)$ 
and $I = \rres(T)$. Moreover, each ideal $\rres_k(T)$ is cofinite and $\tau_{\rres(T)}$ is a profinite congruence on 
$A^{-\omega}$, indeed the intersection of the open congruences $\tau_{\rres_k(T)}$. 

Using the left-right duals of Sections~\ref{section.free profinite}-\ref{section.src omega}, we have similar results 
for $\lres(T)$ and the sequence $(\lres_k(T))$.

In a subsequent paper, we intend to characterize polynomial time Turing machines in this framework, including the natural semaphore codes action and the action of $\beta^{(n)}$ and $\beta^{(\infty)}$. The approach 
will constitute a variation of \cite{RW:1995}: we will need to consider certain metrics that will give the same topology 
as in Sections~\ref{section.free profinite}-\ref{section.src omega}, but conditions involving the metrics will take us 
from the realm of topology into that of geometry.


\end{document}